\documentclass[11pt,a4paper]{article} 
\usepackage[cp1251]{inputenc}         
\usepackage{amssymb,amsfonts,amsmath} 
\usepackage{graphicx}                 
\usepackage{enumerate}                %

\newtheorem{ter}{Theorem}[section]

\newtheorem{pr}{Proposition}[section]
\newtheorem{cor}{Corollary}[section]
\newtheorem{lm}{Lemma}[section]
\newtheorem{rem}{Remark}[section]
\newtheorem{quest}{Question}[section]
\begin{document}

\title{Affine quotients of supergroups}
\author{A.N.Zubkov}
\date{}

\maketitle
\begin{abstract}
In this article we consider sheaf quotients of affine superschemes
by affine supergroups that act on them freely. The necessary and
sufficient conditions for such quotients to be affine are given.
If $G$ is an affine supergroup and $H$ is its normal
supersubgroup, then we prove that a dur $K$-sheaf
$\tilde{\tilde{G/H}}$ is again affine supergroup. Additionally, if
$G$ is algebraic, then a $K$-sheaf $\tilde{G/H}$ is also algebraic
supergroup and it coincides with $\tilde{\tilde{G/H}}$. In
particular, any normal supersubgroup of an affine supergroup is
faithfully exact.
\end{abstract}

\section*{Introduction}

Let $G$ be an algebraic group defined over an algebraically closed
field $K$ of an arbitrary  characteristic. A closed subgroup $H$
of $G$ is called {\it exact} iff the induced functor $ind^G_H$ is
exact. A remarkable theorem of Cline-Parshall-Scott says that $H$
is exact iff the quotient $G/H$ is affine iff $K[G]$ is an
injective $H$-module \cite{cps}. In the category of affine
schemes, well-known Takeuchi's theorem states that the dur
$K$-sheaf (faisceau dur in the terminology from \cite{dg}) of
right cosets $\tilde{\tilde{G/H}}$ is affine iff $K[G]$ is a
faithfully coflat right (or/and left, respectively)
$K[H]$-comodule iff $K[G]$ is an injective cogenerator in the
category of right (or/and left, respectively) $K[H]$-comodules
\cite{t2, t4}. In this case $H$ is called {\it faithfully exact}.

In the category of affine supergroups, only the second equivalence
of the above Cline-Parshall-Scott theorem has been recently proved
in \cite{z}. The definition of a dur $K$-sheaf (as well as the
definition of a $K$-sheaf or faisceau, see \cite{dg}) can be
easily adapted to the category of $K$-functors over commutative
superalgebras. Therefore, it is natural to ask whether, for an
affine supergroup $G$  and its supersubgroup $H$, the dur
$K$-sheaf $\tilde{\tilde{G/H}}$ is affine. In the case when $G$ is
algebraic, that is $K[G]$ is a finitely generated superalgebra,
one can also ask whether the $K$-sheaf $\tilde{G/H}$ is affine.
Notice that in general we only have inclusion
$\tilde{G/H}\subseteq \tilde{\tilde{G/H}}$ but if $\tilde{G/H}$ is
an affine superscheme, then $\tilde{G/H}=\tilde{\tilde{G/H}}$.

In the present article we consider more general case when $G$ acts
on an affine superscheme $X$ freely and we give necessary and
sufficient conditions for $\tilde{\tilde{X/G}}$ (and
$\tilde{X/G}$, respectively) to be an affine superscheme.
Following \cite{mw, t2, t4} we generalize some theorems about Hopf
algebras and their coideal subalgebras to Hopf superalgebras and
their coideal supersubalgebras. In particular, we obtain a
superalgebra version  of the above Takeuchi's theorem.

Equipped with these results, we proceed to prove the main theorem
of our article. Namely, the fundamental theorem of the algebraic
group theory states that if $G$ is an algebraic group and $H$ is
its closed normal subgroup, then $G/H$ is again an algebraic
group, see \cite{dg, ham, w}. For any affine supergroup $G$ and
its closed normal supersubgroup $H$ we prove that
$\tilde{\tilde{G/H}}$ is again an affine supergroup. Moreover, if
$G$ is algebraic, then $\tilde{\tilde{G/H}}=\tilde{G/H}$ is also
algebraic. It follows that any closed normal supersubgroup of an
affine supergroup is faithfully exact.

The article is organized as follows. In the first section
necessary definitions and notations concerning super(co)algebras
and super(co)modules over them are collected. On top of that,
(co)flat and faithfully (co)flat super(co)modules are defined and
some well-known results about flat and faithfully flat modules
over commutative algebras are generalized to (super)modules over
commutative superalgebras. In the second section we introduce
$K$-functors as functors from the category of commutative
superalgebras to the category of sets. More attention is paid to
the particular case of affine superschemes. We also define (dur)
$K$-sheafs and completions of certain $K$-functors in the
Grothendieck topology of faithfully flat (faithfully flat and
finitely presented) coverings. In the third section we define the
superspace of distributions of an affine superscheme. This
superspace has a natural structure of (cocommutative) Hopf
superalgebra provided the original superscheme is a supergroup.
Moreover, if $char K=0$ and this supergroup is algebraic, then we
prove that the corresponding distribution superalgebra is
isomorphic to the universal enveloping superalgebra of its Lie
superalgebra. In addition, we introduce the notion of {\it Lie
superalgebra functor} by means of superalgebra of dual numbers
(see \cite{f1}).

In the fourth section we consider an affine superscheme $X$ and a
supergroup $G$ that acts on $X$ (on the right) freely. We give
necessary and sufficient conditions for a quotient
$\tilde{\tilde{X/G}}$ ($\tilde{X/G}$) to be affine. The main
results of the fifth section are Theorem 5.1 and 5.2. Theorem 5.1
says surprisingly more about Hopf supersubalgebras than one would
expect. Keeping in mind Takeuchi's theorem it would be quite
natural to expect that a Hopf superalgebra is (left and right)
faithfully flat supermodule over its Hopf supersubalgebra but it
is actually faithfully flat as a module!

In sixth and ninth sections we prove the main result of this
article. The principal idea is different from \cite{t1, ham} for
the following reason. In the category of vector superspaces, an
exterior power does not have the same nice properties as in the
category of vector spaces. More precisely, if $V$ is a superspace
of finite (super)dimension $(m, n)$ and $W$ is its supersubspace
of (super)dimension $(s, t)$, where $t > 0$, then no exterior
power $\Lambda^d(W)$ is one-dimensional. In particular, the
"naive" or "direct" generalization of Chevalley's theorem is not
possible. Instead of exterior powers one can use a {\it
berezinian} $Ber_{s|t}(W)$ of $W$ but it does not solve our
problem. In fact, we would have to embed $Ber_{s, t}(W)$ to
something like $Ber_{s, t}(V)$ (analogously as $\Lambda^k(U)$ is
embedded into $\Lambda^k(L)$, where $U$ is a subspace of a space
$L$ and $k=\dim U$) but there is no appropriate definition of such
a supermodule. To overcome this obstacle, we construct to any
normal supersubgroup some biggest normal supersubgroup over which
the quotient is affine. Next step is to prove that the original
supersubgroup coincides with this new one. For $char K =p > 0$ we
use some trick with a Frobenius map. The characteristic zero case
is much more complicated and we have to introduce a notion of a
{\it pseudoconnected component} of a supergroup to reduce our
problem to the case of a finite normal supersubgroup. Besides, we
use induction on superdimensions of Lie superalgebras of our
supergroups and some properties of adjoint representations. In the
final section an example of faithfully exact supersubgroup is
given. This is a Levi supersubgroup of a general linear
supergroup. Furthermore, we construct an affine superscheme $X$ on
which a finite (odd) unipotent supergroup $G$ acts in a such way,
that neither $\tilde{\tilde{X/G}}$, nor $\tilde{X/G}$ is affine.

In seventh and eighth section a partial answer for the following
Brundan's question is given. Let $G$ be an algebraic supergroups
and $H$ be its supersubgroup such that $H_{ev}$ is reductive. Is
it true that $\tilde{\tilde{G/H}}$ is affine? We show that the
answer is positive if $char K=p >0$ or $G$ is finite.

\section{Super(co)algebras and super(co)modules}

We follow definitions and notations from \cite{z, mu} (see also
\cite{l}). Let $K$ be a field of characteristic $p\neq 2$. For a
$K$-vector superspace $V$ its {\it superdimension} is defined as
$\mbox{s}\dim V= (\dim V_0, \dim V_1)$. Let $A$ be a (associative)
superalgebra over $K$. Denote by $A-smod$ (or $smod-A$,
respectively) the category of all left (or right, respectively)
$A$-supermodules with even morphisms. Let $X$ be a left (or right)
$A$-supermodule with ${\bf Z}_2$-grading $X=X_0\oplus X_1$. For
any $x\in X$ denote by $x_0$ and $x_1$ its homogeneous components,
that is $x_0\in X_0, x_1\in X_1$ and $x=x_0+x_1$. If $X\in smod-A,
Y\in A-smod,$ then the tensor product $X\otimes _A Y$ has a
natural ${\bf Z}_2$-grading given by $|x\otimes y|= |x|+|y|\pmod
2$ for $x\in X, y\in Y$.

A superalgebra $A$ has the opposite companion $A^{\circ}$ whose
underlying superspace coincides with $A$ and the multiplication in
$A^{\circ}$ is defined by $a\ast b=(-1)^{|a||b|}ba$. It is clear
that $(A^{\circ})^{\circ}=A$. We have an equivalence of categories
$A-smod\simeq smod-A^{\circ}$ given by $M\mapsto M^{\circ}$ for
$M\in A-smod$, where $M^{\circ}$ coincides with $M$ as a
superspace and a structure of a right supermodule on $M^{\circ}$
is defined by $m\ast a=(-1)^{|m||a|}am$ for $a\in A, m\in M$. For
any $X\in smod-A, Y\in A-smod,$ we have an isomorphism $X\otimes_A
Y\simeq Y^{\circ}\otimes_{A^{\circ}} X^{\circ}$ given by $x\otimes
y\mapsto (-1)^{|x||y|}y\otimes x$ for $x\in X, y\in Y$ which is
functorial in $X$ and $Y$.

A superalgebra $A$ is called {\it commutative} if any homogeneous
elements $a, b\in A$ satisfy $ab=(-1)^{|a||b|}ba$. In particular,
$A=A^{\circ}$ and $A-smod\simeq smod-A$. In other words, any
one-sided $A$-supermodule has a canonical structure of a
$A$-superbimodule. Denote by $SAlg_K$ the category of all
commutative $K$-superalgebras with even morphisms. If $\phi : B\to
A, \psi : B\to C$ are morphisms in $SAlg_K$, then $A\otimes _B
C\in SAlg_K$ (see \cite{l}). From now on, all superalgebras are
commutative unless otherwise stated.

Let $K[m|n]=K[t_1 ,\ldots , t_m | z_1,\ldots z_n]$ be a free
commutative superalgebra with free generators $t_1 ,\ldots , t_m ,
z_1,\ldots z_n$, where $ |t_i|=0, |z_j|=1$ for  $1\leq i\leq m,
1\leq j\leq n$. It can be identified with the symmetric
superalgebra $S(V)$ of a superspace $V$, where $\dim V_0 =m, \dim
V_1 =n$ (see \cite{z}).
\begin{lm} Let $R=R_0\bigoplus R_1$ be a superalgebra. Then \\
i) If $R$ is finitely generated, then $R$ is noetherian; \\
ii) Every one-sided superideal of $R$ is two-sided; \\
iii) A prime ideal $\cal P$ of $R$ has a form ${\cal P}={\cal P}_0
+ R_1$, where ${\cal P}_0$ is a prime ideal of $R_0$. If $\cal M$
is a maximal one-sided ideal of $R$, then ${\cal M}={\cal M}_0
+R_1$, where ${\cal M}_0$ is a maximal ideal of $R_0$. In
particular, all these ideals are superideals.
\end{lm}
Proof. The statement ii) is obvious. To prove i) we notice that
$K[m|n]$ is finitely generated module over $K[t_1 ,\ldots , t_m]$.
Finally, iii) holds since $RR_1$ is a nil ideal.
\begin{rem}
It is not true that every one-sided ideal in a superalgebra is
necessary two-sided. For example, the left ideal $K[1|2](t_1+z_1)$
is not two-sided.
\end{rem}
Denote by  $\sqrt[R]{I}$ the prime radical of a superideal $I$. It
is the intersection of all prime ideals containing $I$. It is not
difficult to see that
$$\sqrt[R]{I}=\{r\in R | \exists n, r^n\in I\}=\{r\in R | \exists n, r_0^n\in I_0\},$$
where for the last equality we used the obvious formula $r^n=r_0^n
+nr_0^{n-1}r_1$.

Let $A$ be a (not necessary commutative) superalgebra. A left
supermodule $Y\in A-smod$ (a right supermodule $Y\in smod-A$,
respectively) is called {\it flat} if the functor $X\to X\otimes
_A Y$ ($X\to Y\otimes_A X$, respectively) is an exact functor from
the category $smod-A$ ($A-smod$, respectively) to the category of
superspaces.

Next, $Y\in A-smod$ ($Y\in smod-A$, respectively) is called {\it
faithfully flat} if the corresponding functor is faithfully exact,
that is the exactness of any sequence $X'\to X\to X''$ in $smod-A$
(in $A-smod$, respectively) is equivalent to the exactness of the
sequence of superspaces $X'\otimes_A Y\to X\otimes_A Y\to
X''\otimes_A Y$ ($Y\otimes_A X'\to Y\otimes_A X\to Y\otimes_A
X''$, respectively). Using the previous discussion, a left
$A$-supermodule $Y$ is flat (faithfully flat) iff the right
$A^{\circ}$-supermodule $Y^{\circ}$ is flat (faithfully flat). If
$A$ is commutative, then any supermodule is left flat (left
faithfully flat) iff it is right flat (right faithfully flat).
Most of standard characterizations of flatness or faithful
flatness from \cite{bur} can be easily translated to the
supercase. We call such translation a {\it superversion} of the
corresponding statement. Proofs of superversions of results from
\cite{bur} that are not difficult are left to the reader.

Let $A$ be an algebra and let $S$ be a multiplicative set
belonging to the center of $A$. The algebra of fractions $S^{-1}A$
and the left (or right, respectively) $S^{-1}A$-module of
fractions $S^{-1}X\simeq S^{-1}A\otimes _A X$ (or $S^{-1}X\simeq
X\otimes_A S^{-1}A$, respectively) for left (or right,
respectively) $A$-module $X$ is defined in the usual way. If $A$
is a superalgebra, $X$ is an $A$-supermodule and $S\subseteq A_0$,
then $S^{-1}X$ is also an $A$-supermodule with ${\bf Z}_2$-grading
given by $(S^{-1}X)_i= S^{-1}X_i$ for  $i=0, 1$.
\begin{lm} Let $A$ and $S$ be as above. Then \\
i) $S^{-1}A$ is a flat $A$-module; \\
ii) If central elements $a_1,\ldots , a_n\in A$ generate $A$ as an
ideal, then the algebra $\prod_{1\leq i\leq n}A_{a_i}$ is a
faithfully flat (left and right) $A$-module.
\end{lm}
Proof. The first statement is an easy generalization of Theorem 1
from \cite{bur}, II, \S 2. To prove the second statement, use
Proposition 1 from \cite{bur}, I, \S 3 and observe that all powers
of the elements $a_1,\ldots , a_n$ again generate $A$ as an ideal.
\begin{lm}
A superalgebra $A$ is generated by elements $a_1,\ldots , a_n$ as
a left (or right) ideal iff $A$ is generated by their even
components.
\end{lm}
Proof. Assume that $1= \sum_{1\leq i\leq n} b_i a_i$. Set $a_{i,
k}, b_{i, k}\in A_k ,$ for $ 1\leq i\leq n$ and $k=0, 1$. Then
$$\sum_{1\leq i\leq n}b_{i, 0}a_{i, 0}=1-
\sum_{1\leq i\leq n}b_{i, 1}a_{i, 1}$$ and the element
$1-\sum_{1\leq i\leq n}b_{i, 1}a_{i, 1}\in 1+AA_1$ is invertible.
\begin{cor}
If $A$ is a superalgebra and $a_1,\ldots , a_n\in A_0$ generate
$A_0$ as ideal, then the superalgebra $\prod_{1\leq i\leq
n}A_{a_i}$ is a faithfully flat (left and right) $A$-module.
\end{cor}
The spectrum of all maximal ideals of a superalgebra $A$ is
denoted by $Max(A)$. For any ${\cal M}\in Max(A)$ we denote by
$N_{{\cal M}}$ an {\it even} localization of an $A$-supermodule
$N$. More precisely, $N_{{\cal M}}=(A_0\setminus {\cal M}_0)^{-1}
N$.

In what follows all algebras are superalgebras.
\begin{lm}
If ${\cal M}\in Max(A)$, then  the algebra $A_{{\cal M}}$ is local
and ${\cal M}A_{{\cal M}}$ is its Jacobson radical. In particular,
the left (and right) $A$-module $B=\bigoplus_{{\cal M}\in Max(A)}
A_{{\cal M}}$ is faithfully flat.
\end{lm}
Proof. If an element $a=a_0+ a_1\in A$ is such that $a_0\not\in
{\cal M}_0$, then $a^{-1}=\frac{1}{a_0}-\frac{a_1}{a_0^2}\in
A_{{\cal M}}$. Next, by Lemma 1.2, the $A$-module $B$ is flat and
${\cal M}B\neq B, B{\cal M}\neq B$ for all ${\cal M}\in Max(A)$.
\begin{lm}
A morphism of left (or right) $A$-modules $M\to N$ is a
monomorphism (an epimorphism or an isomorphism, respectively) iff
the induced morphism $M_{{\cal M}}\to N_{{\cal M}}$ is injective
(surjective or bijective, respectively) for any ${\cal M}\in
Max(A)$.
\end{lm}
Proof. A word-by-word repetition of the proof of Theorem 1, II, \S
3, \cite{bur} combined with Lemma 1.4.
\begin{lm}
A left (or right) $A$-module $M$ is flat (faithfully flat,
respectively) iff $M_{{\cal M}}$ is a flat (faithfully flat,
respectively) $A_{{\cal M}}$-module for all ${\cal M}\in Max(A)$.
\end{lm}
Proof. The necessary condition is a consequence of Proposition 8,
I, \S 2 and Proposition 4, I, \S 3 from \cite{bur}. For the
sufficient condition, assume that $M_{{\cal M}}$ is a flat (and,
for example, right) $A$-module for all ${\cal M}\in Max(A)$. If
$N_1\to N_2$ is an inclusion of left $A$-modules, combine
$(M\otimes_A N)_{{\cal M}}\simeq M_{{\cal M}}\otimes_A N$ with
Lemma 1.5 to obtain that $M\otimes_A N_1\to M\otimes_A N_2$ is
again an inclusion. Additionally, if $M_{{\cal M}}$ is faithfully
flat for all ${\cal M}\in Max(A)$, then ${\cal M} M_{{\cal M}}\neq
M_{{\cal M}}$ implies ${\cal M}M\neq M$.
\begin{pr}
Let $\phi : A\to B$ be a morphism of superalgebras and let $M$ be
a $B-A$-bimodule (an $A-B$-bimodule, respectively) such that
$ma=\phi(a)m$ for $m\in M, a\in A_0$ ($am=m\phi(a)$ for $m\in M,
a\in A_0$, respectively). Then the
following properties are equivalent: \\
i) $M$ is a flat $A$-module; \\
ii) $M_{{\cal N}}$ is a flat $A$-module for all ${\cal N}\in Max(B)$; \\
iii) $M_{{\cal N}}$ is a flat $A_{{\cal M}}$-module for every
${\cal N}\in Max(B)$, where ${\cal M}=\phi^{-1}({\cal N})$.
\end{pr}
Proof. Use Lemmas 1.2 and 1.5 to copy the proof of Proposition 15,
II, \S 3, \cite{bur}.
\begin{cor}
If $M$ is a $B$-supermodule, then the conditions of Proposition
1.1 hold automatically. In particular, the properties i)-iii) are
equivalent.
\end{cor}

Let $M$ be a flat left (or right, respectively) $A$-module over an
algebra $A$ and let $\alpha$ be an automorphism of $A$. Denote by
$M^{\alpha}$ an $A$-module such that $M^{\alpha}=M$ and $a\star
m=\alpha(a)m $ (or $m\star a=m\alpha(a)$, respectively) for $a\in
A, m\in M$. The corollary after Proposition 13, I, \S 2 of
\cite{bur} implies that $M$ is a flat $A$-module iff $M^{\alpha}$
is a flat $A$-module.

All necessary definitions and notations concerning
supercoalgebras, Hopf superalgebras and supercomodules over them
can be found in \cite{z, mu}. If $C$ is a supercoalgebra  and $V$
is a left (or right, respectively) $C$-supercomodule, then a
counit of $C$ and a coaction map $V\to C\otimes V$ (or $V\to
V\otimes C$, respectively) are denoted by $\epsilon_C$ and
$\tau_V$ correspondingly. A comultiplication of $C$ is denoted by
$\delta_C$. Additionally, if $C$ is a Hopf superalgebra, then its
antipode is denoted by $s_C$. The category of left (or right,
respectively) $C$-supercomodules with even morphisms is denoted by
$C-scomod$ (or $scomod-C$, respectively). If $V\in scomod-C, W\in
C-scomod$, then one can define a {\it cotensor} product
$$V\square_C W=\{x\in V\otimes W | (\tau_V\otimes
id_W)(x)=(id_V\otimes\tau_W)(x)\}.$$ A left (or right,
respectively) $C$-supercomodule $V$ is called {\it (faithfully)
coflat} if the functor $W\to W\square_C V$ (or $W\to V\square_C
W$, respectively) is (faithfully) exact, where $W\in scomod-C$ (or
$W\in C-scomod$, respectively).
\begin{lm} (see A.2.1, \cite{t4}) A right (or left) $C$-supercomodule $V$ is coflat
(faithfully coflat, respectively) iff $V$ is injective (an
injective cogenerator, respectively).
\end{lm}
Proof. Let $V\in scomod-C$, $W\in C-scomod$ and $\dim W  <
\infty$. The dual superspace $W^*$ has a uniquely defined
structure of a right $C$-supercomodule given by $\sum
f_1(w)c'_2=\sum f(w_1)c_2,$ where $\tau_W(w)=\sum c_2\otimes w_1,
\tau_{W^*}(f)=\sum f_1\otimes c'_2$ for $w\in W, f\in W^*$. We
have an isomorphism of superspaces $V\square_C W\to Hom_C(W^* ,
V)$ induced by $v\otimes w(f)=f(w)v$ for $v\in V$, $w\in W, f\in
W^*$. Since this isomorphism is functorial in $W$, we conclude the
proof as in \cite{t4}.

\section{$K$-functors and $K$-sheafs (faisceaux)}

Following the book \cite{jan} we call a functor from the category
$SAlg_K$ to the category of sets a $K$-functor. The category of
all $K$-functors is denoted by $\cal F$. A $K$-functor $SSp \ R$
defined as $SSp \ R(A)=Hom_{SAlg_K}(R, A)$ for $A\in SAlg_K$ is
called an {\it affine superscheme} (this definition is different
from the definition used in \cite{z} since we do not suppose that
$R$ is finitely generated). The superalgebra $R\in SAlg_K$ is
called a {\it coordinate superalgebra} of the superscheme $SSp \
R$. If $X=SSp \ R$, then $R$ is also denoted by $K[X]$.
\begin{lm}(Yoneda's lemma, \cite{jan}, part I, (1.3))
For an affine superscheme $SSp \ R$ and a $K$-functor $X$ there is
a canonical isomorphism $Mor(SSp \ R , X)\simeq X(R)$ which is
functorial in both arguments. In particular, the category $SAlg_K$
is anti-equivalent to the full subcategory of affine superschemes.
\end{lm}
Proof. The statement of the lemma is a partial case of more
general theorem about covariant representable functors (see
\cite{bd}, Theorem 1.6). The isomorphism is given by $f\mapsto
x_f=f(R)(id_R)$ for $f\in Mor(SSp \ R , X)$ and the inverse map is
given by $x\mapsto f_x$, where $f_x(\alpha)=X(\alpha)(x)$ for
$x\in X(R), \alpha\in Hom_{SAlg_K}(R, A)$ and $A\in SAlg_K$.
\begin{cor}
The universal property of $A\otimes_B C$ implies a canonical
isomorphism $$SSp \ A\times_{SSp \ B} SSp \ C\simeq SSp \
A\otimes_B C .$$
\end{cor}
The affine superscheme ${\bf A}^{m|n}=SSp \ K[t_1 ,\ldots , t_m |
z_1,\ldots z_n]$ is called {\it $(m|n)$-affine superspace}. It is
clear that ${\bf A}^{m|n}(B)=B_0^m\bigoplus B_1^n$ for $B\in
SAlg_K$. In particular, ${\bf A}^{1|1}(B)=B$ for any superalgebra
$B$.

Let $I$ be a superideal of $R\in SAlg_K$. Denote by $V(I)$ a {\it
closed subfunctor} of $SSp \ R$ corresponding to $I$. By
definition, $V(I)(A)=\{\phi\in SSp \ R(A)|\phi(I)=0\}$. It is
obvious that $V(I)\simeq SSp \ R/I$. All standard properties of
closed subfunctors of affine schemes mentioned in \cite{jan}, part
I, (1.4) are translated to the category of affine superschemes per
verbatim.

Let $X$ be an affine superscheme. A functor $Y\subseteq X$ is
called {\it open} if
$$Y(A)=\{x\in X(A)|\sum_{f\in I}Ax(f)=A\}$$
for a subset $I\subseteq K[X]$ and $A\in SAlg_K$. Denote this
functor by $D(I)$.
\begin{lm}
i) If $J$ is the smallest superideal
containing $I$, then $D(I)=D(J)=D(J_0)$; \\
ii) Let $I$ and $I'$ be superideals of $R$. Then $D(I)\subseteq
D(I')$ iff $\sqrt[R]{I}\subseteq \sqrt[R]{I'}$. Additionally,
$D(I)=D(\sqrt[R]{I})=D(\sqrt[R_0]{I_0})$.
\end{lm}
Proof. All statements can be proved by the same trick with a
representation of unit as in Lemma 1.3 and by the standard
reductions to quotients modulo prime ideals (see \cite{jan}, part
I (1.6)).

An important example of an open subfunctor is a so-called {\it
principal open} subfuctor $X_f=D(\{f\})$ for $f\in K[X]$. It can
be checked easily that $X_f=X_{f_0}$ is again an affine
superscheme and $K[X_f]=K[X]_{f_0}$. All other properties of open
subfunctors mentioned in \cite{jan}, part I (1.6) are easily
translated to the category of affine superschemes.

Let $G$ be a {\it group} $K$-functor, that is a $K$-functor to the
category of groups. We say that $G$ acts on a $K$-functor $X$ on
the right (on the left, respectively) if there is a morphism of
functors $f : X\times G\to X$ ($g : G\times X\to X$, respectively)
such that $f (id_X\times\mu)=f(f\times id_G)$ and $f i_E=id_X$
($g(\mu\times id_X)=g(id_G\times g)$ and $g j_E =id_X$,
respectively). Here $\mu : G\times G\to G$ is a multiplication of
$G$ and $i_X : X\to X\times G$ is defined as $i_X(R)(x)=(x,
1_{G(R)})$ ($j_X(R)(x)=(1_{G(R)}, x)$, respectively) for $x\in
X(R), R\in SAlg_K$. From now on we consider any action on right
unless otherwise stated.

It is obvious that the category of affine group superschemes (=
affine supergroups) is anti-equivalent to the category of
commutative Hopf superalgebras (see \cite{z, jan}). If $G$ is an
affine supergroup, then denote by $\epsilon_G, \delta_G$ and $s_G$
the counit, comultiplication and antipode of $K[G]$
correspondingly. If $K[G]$ is finitely generated, then $G$ is
called an {\it algebraic} supergroup. Closed supersubgroups $H\leq
G$ are in one-to-one correspondence with Hopf superideals
$I_H\subseteq K[G]$ such that $H=V(I_H)$. If $X$ and $G$ are
affine, then an action of $G$ on $X$ is uniquely defined by a
morphism of superalgebras $\tau : K[X]\to K[X]\otimes K[G]$ with
respect to which $K[X]$ is a right $K[G]$-supercomodule.

Let $V$ be a vector superspace of superdimension $(m, n)$. Denote
by $GL(V)$ or by $GL(m|n)$ the corresponding general linear
supergroup. More precisely, $GL(V)$ is a group $K$-functor such
that for any $B\in SAlg_K$ the group $GL(V)(B)$ consists of all
even and $B$-linear automorphisms of $V\otimes B$. It is not
difficult to see that $GL(V)$ is an algebraic supergroup (see
\cite{z, bk, bkl} for more details). A {\it linear representation}
of a group $K$-functor $G$ is a morphism of group $K$-functors
$\rho :
 G\to GL(V)$. In this case the superspace $V$ is called $G$-{\it supermodule}.
If $G$ is an affine supergroup, then $V$ is a $G$-supermodule iff
it is a right $K[G]$-supercomodule (see \cite{z, jan}). In fact,
fix a basis $v_1,\ldots , v_{m+n}$ of $V$ such that $|v_i|=0$ if
$1\leq i\leq m$ and $|v_i|=1$ otherwise. Set
$\tau_V(v_i)=\sum_{1\leq j\leq m+n} v_j\otimes r_{ji}$ for $1\leq
i\leq m+n$. Then $\rho(g)=(g(r_{ij}))$ for $g\in G(B)$ and $B\in
SAlg_K$.

From now on, any group $K$-functor is affine unless otherwise
stated. Let $W$ be a supersubspace of a finitely dimensional
$G$-supermodule $V$. The stabilizer $Stab_G(W)$ is a group
subfunctor defined as $Stab_G(W)(A)=\{g\in G(A) | g(W\otimes
1)\subseteq W\otimes A\}$ for $A\in SAlg_K$. It is easy to see
that $Stab_G(W)$ is a closed supersubgroup of $G$. In fact,
without loss of generality one can assume that $v_1 ,\ldots, v_s,
v_{m+1},\ldots , v_{m+t}$ for $s\leq m, t\leq n$ is a basis of
$W$. Denote by $M$ the set of indexes $\{1,\ldots, s, m+1,\ldots
m+t\}$. Then $g\in Stab_G(W)(A)$ iff $g(r_{ji})=0$ for all
$j\not\in M, i\in M$ and $A\in SAlg_K$.

Let $R_1, \ldots , R_n$ be a finite family of commutative
$R$-superalgebras with respect to a set of morphisms
$\iota^{R_i}_R : R\to R_i$ in $SAlg_K$. Such a family is called
{\it faithfully flat covering} of $R$ (ff-covering, for short)
whenever $R$-supermodule $R_1\times\ldots\times R_n$ is faithfully
flat. We say that $R$-superalgebra $R'$ is {\it finitely
presented} if $R'\simeq R[m|n]/I$, where $R[m|n]=R\otimes K[m|n]$
and $I\subseteq R[m|n]$ is a finitely generated superideal. It is
not difficult to check that $R'$ is a finitely presented
$R$-superalgebra iff $R'\simeq R\otimes_{A} A[m|n]/I$, where $A$
is a finitely generated supersubalgebra of $R$ (see \cite{dg}, I,
\S 3). Following \cite{jan} we call a ff-covering $R_1,\ldots ,
R_n$ {\it fppf-covering} if all $R_i$ are finitely presented
$R$-superalgebras.

A $K$-functor $X$ is called {\it dur $K$-sheaf} (or {\it faisceau
dur}) if for any ff-covering $R_1 ,\ldots R_n$ of a superalgebra
$R$ the diagram
$$X(R)\to\prod_{1\leq i\leq n}X(R_{i})\begin{array}{c}\to \\
\to\end{array}\prod_{1\leq i, j\leq n} X(R_i\otimes_R R_j)$$ is
exact, where the last two maps are induced by morphisms $R_i\to
R_i\otimes_R R_j$ and $R_i\to R_j\otimes_R R_i$, respectively,
defined as $a\mapsto a\otimes 1$ and $b\mapsto 1\otimes b$ for $a,
b\in R_i$. This property is equivalent to the following two
conditions. For all $R_1 ,\ldots , R_n, R, R'\in SAlg_K$, where
$R'$ is a faithfully flat $R$-supermodule, there is a canonical
bijection $X(\prod_{1\leq i\leq n}R_i)\simeq\prod_{1\leq i\leq
n}X(R_i)$ and the diagram
$$X(R)\to X(R')\begin{array}{c}\to \\
\to\end{array}X(R'\otimes_R R')$$ is exact, see \cite{jan}.
Replacing ff-coverings by fppf-coverings we obtain a definition of
a {\it $K$-sheaf} (or {\it faisceau}), cf. \cite{dg, jan}. Denote
the full subcategory of $K$-sheafs (dur $K$-sheafs, respectively)
by $\tilde{\cal F}$ ($\tilde{\tilde{\cal F}}$, respectively). It
is clear that $\tilde{\tilde{\cal F}}\subseteq \tilde{\cal F}$ and
it can be checked easily that any affine supersheme is a dur
$K$-sheaf, see \cite{jan}, part I (5.3).

For a $K$-functor $X$ one can construct an {\it associated
$K$-sheaf} $\tilde X$ and a dur $K$-sheaf $\tilde{\tilde{X}}$
following the way described in \cite{dg}, III, \S 1 . We consider
a partial case following \cite{jan}, part I (5.4).
Assume that the $K$-functor $X$ satisfies the following conditions : \\
(*) $X(\prod_{1\leq i\leq n}R_i)\simeq \prod_{1\leq i\leq
n}X(R_i)$ for all $R_1 ,\ldots, R_n\in SAlg_K$; \\
(**) $X(R)\to X(R')$ is an inclusion for arbitrary fppf-covering
$R'$ of superalgebra $R$.

The family of all $K$-functors satisfying the properties (*) and
(**) is closed under direct products. Define a partial order on
the set of all fppf-coverings of a superalgebra $R$ by the
following rule: $R'\leq R''$ if $R''$ is a fppf-covering of $R'$.
Clearly, this poset is directed since $R' , R''\leq R'\otimes_R
R''$ for any two fppf-coverings of $R$. For each superalgebra $R$
define the direct spectrum
$${\cal X}(R)=\{X(R' , R)=\ker(X(R')\begin{array}{c}\to \\
\to\end{array}X(R'\otimes_R R'))| R' \ \mbox{is a fppf-covering
of} \ R\}$$ with canonical inclusions $X(R', R)\to X(R'' , R)$ for
any couple $R'\leq R''$. Observe that $R''\otimes_R R''$ is a
fppf-covering of $R'\otimes_R R'$ and set
$\tilde{X}(R)=\lim\limits_{\rightarrow}{\cal X}(R)$. The functor
$\tilde X$ is the required  completion of $X$ with respect to the
Grothendieck topology of fppf-coverings. A canonical inclusion
$\alpha_X : X\to\tilde X$ induces a canonical bijection
$Mor(\tilde{X}, Y)\to Mor(X, Y)$ for any $K$-sheaf $Y$. One gets
easily that if $X$ is a subfunctor of a $K$-sheaf $Y$ satisfying
(*), then $\tilde X\subseteq Y$. Moreover,
$$\tilde{X}(R)=\{y\in Y(R)|\mbox{there is} \ R'\geq R \
\mbox{such that} \ Y(\iota_R^{R'})(y)\in X(R')\}.$$ To prove all
of the above statements we only need a superversion of Proposition
4 \cite{bur}, I, \S 3 together with an additional statement which
says that $R=\prod_{1\leq i\leq n}R_i\leq R'$ iff $R_i=Re_i\leq R'
e_i$ for any $i$, where $e_i=(0,\ldots ,\underbrace{1}_{i-\mbox{th
\ place}},\ldots , 0)$ and $1\leq i\leq n$. Arguments for dur
$K$-sheafs are the same, except that fppf-coverings are replaced
by ff-coverings.
\begin{lm}
Let $G$ be a group $K$-functor that satisfies (*) and (**) for all
ff-coverings (fppf-coverings, respectively). Then
$\tilde{\tilde{G}}$ ($\tilde{G}$, respectively) is again a group
dur $K$-sheaf (a group $K$-sheaf, respectively) and the canonical
inclusion $G\to \tilde{\tilde{G}}$ ($G\to\tilde{G}$, respectively)
is a morphism of group functors.
\end{lm}
Proof. Let $g, g_1, g_2\in\tilde{\tilde{G}}(A)$. For a suitable
ff-covering $B$ of a superalgebra $A$ we set
$\overline{g}=\tilde{\tilde{G}}(\iota^B_A)(g)$ and $\overline{g}_i
=\tilde{\tilde{G}}(\iota^B_A)(g_i)\in G(B)$ for $i=1, 2$. The pair
$(\overline{g}_1 , \overline{g}_2)$ belongs to $(G\times G)(B, A)$
which implies $\overline{g}_1\overline{g}_2\in G(B, A)$. Define
$g_1
g_2=\tilde{\tilde{G}}(\iota^B_A)^{-1}(\overline{g}_1\overline{g}_2)\in\tilde{\tilde{G}}(A)$
and $g^{-1}=\tilde{\tilde{G}}(\iota^B_A)^{-1}(\overline{g}^{-1})$.
These definitions do not depend on the choice of $B$. If $\phi :
A\to C$ is a morphism in $SAlg_K$, $B$ and $B'$ are ff-coverings
of $A$ and $C$, respectively, then $B\otimes_A B'$ is a
ff-covering of $C$ and $B'$, respectively. Denote by $\alpha$ and
$\beta$ the morphisms $B\to B\otimes_A B'$ and $B'\to B\otimes_A
B'$ given by $\alpha(b)=b\otimes 1, \ \beta(b')=1\otimes b'$ for
$b\in B, b'\in B'$ and set $x_i=\tilde{\tilde{G}}(\phi)(g_i), \
\overline{x}_i=\tilde{\tilde{G}}(\iota^{B'}_C)(x_i)$ for $i=1, 2$.
Then
$$\tilde{\tilde{G}}(\beta\iota^{B'}_C\phi)(g_1 g_2)=\tilde{\tilde{G}}(\alpha\iota^B_A)(g_1
g_2)=G(\alpha)(\overline{g}_1)G(\alpha)(\overline{g}_2)=G(\beta)(\overline{x}_1)G(\beta)(\overline{x}_2)=
\tilde{\tilde{G}}(\beta\iota^{B'}_C)(x_1 x_2)$$ and since the map
$\tilde{\tilde{G}}(\beta\iota^{B'}_C)$ is injective we infer that
$\tilde{\tilde{G}}(\phi)(g_1 g_2)=x_1 x_2$. The remaining
statements of the lemma are now obvious. The case of $K$-sheafs is
anologous.

\section{Superalgebras of distributions and Lie superalgebras}

Let $X$ be an affine superscheme. Following \cite{jan} we call any
element of $Dist_n(X, {\cal M})=(K[X]/{\cal M}^{n+1})^*$ a {\it
distribution} on $X$ with support at ${\cal M}\in Max(K[X])$ of
orded $\leq n$ (notice that $Max(K[X])$ is obviously identified
with $X(K)$). We have $\bigcup_{n\geq 0}Dist_n(X, {\cal
M})=Dist(X, {\cal M})\subseteq K[X]^*$. If $g : X\to Y$ is a
morphism of affine superschemes, then it induces a morphism of
filtered superspaces $dg_{{\cal M}} : Dist(X, {\cal M})\to Dist(Y,
(g^*)^{-1}({\cal M}))$. In particular, if $X=V(I)$ is a closed
supersubscheme of $Y$, then $Dist(X, {\cal M})$ is identified with
a filtered supersubspace $\{\phi\in Dist(Y, {\cal M}) |
\phi(I)=0\}$, where $I\subseteq {\cal M}$.

If $X$ is an algebraic supergroup and ${\cal M}=\ker\epsilon_X$,
then $Dist(X, {\cal M})$ is denoted by $Dist(X)$. In this case
$Dist(X)$ has a structure of a Hopf superalgebra with a
multiplication
$\phi\psi(f)=\sum(-1)^{|\psi||f_1|}\phi(f_1)\psi(f_2)$ for $\phi,
\psi\in Dist(X), f\in K[X]$ and $\delta_X(f)=\sum f_1\otimes f_2,$
with a unit $\epsilon_X$ and with a counit $\epsilon_{Dist(X)} :
\phi\mapsto\phi (1)$. The comultiplication of $Dist(X)$ is dual to
the multiplication of $K[X]$, cf. \cite{bk, bkl}. Finally, an
antipode $s_{Dist(X)}$ is defined by
$s_{Dist(X)}(\phi)(f)=\phi(s_X(f))$ for $\phi\in Dist(X)$ and
$f\in K[X]$.

We have $Dist_k(X)Dist_l(X)\subseteq Dist_{k+l}(X)$ for all $k,
l\geq 0,$ that is the superalgebra $Dist(X)$ is a filtered
algebra. The superspace $Lie(X)=\{\phi\in Dist_1(X) |\phi(1)=0\}$
has a Lie superalgebra structure given by $[\phi , \psi]=\phi\psi
-(-1)^{|\phi||\psi|}\psi\phi$. As a Hopf superalgebra, $Dist(X)$
is cocommutative which means that
$\delta_{Dist(X)}(\phi)=\sum\phi_1\otimes\phi_2=\sum(-1)^{|\phi_1||\phi_2|}\phi_2\otimes\phi_1$.
Additionally, each $Dist_n(X)$ is a supersubcoalgebra of
$Dist(X)$. For arbitrary morphism of algebraic supergroups $g :
X\to Y$ its {\it differential} $dg=dg_{{\cal M}} : Dist(X)\to
Dist(Y)$ is a homomorphism of filtered Hopf superalgebras. In
particular, its restriction to $Lie(X)$ is a homomorphism of Lie
superalgebras.

Let $L$ be a Lie superalgebra. Denote by $U(L)$ its {\it
(universal) enveloping} superalgebra, see \cite{sch}. The
superalgebra $U(L)$ is a Hopf superalgebra with a comultiplication
defined by $\delta_{U(L)}(x)=x\otimes 1 +1\otimes x$ for $x\in L$.
Its counit $\epsilon_{U(L)}$ is defined by $\epsilon_{U(L)}(L)=0$
and its antipode $s_{U(L)}$ is defined by $s_{U(L)}(x)=-x$ for
$x\in L$. The antipode $s_{U(L)}$ is an anti-automorphism of
$U(L)$ such that
$s_{U(L)}(uv)=(-1)^{|u||v|}s_{U(L)}(v)s_{U(L)}(u)$ for $u, v\in
U(L)$. The Hopf superalgebra $U(L)$ is obviously filtered and
cocommutative.

The inclusion $Lie(X)\subseteq Dist(X)$ induces a morphism $g_X :
U(Lie(X))\to Dist(X)$ of superalgebras.
\begin{lm}(see Lemma 1.2, II, \S 6, \cite{dg})
If $char K =0$ and $X$ is an algebraic supergroup, then $g_X$ is
an isomorphism of Hopf superalgebras.
\end{lm}
Proof. Let $\phi_1 ,\ldots ,\phi_n$ form a basis of $Lie(X)$ dual
to a homogeneous basis $f_1 ,\ldots , f_n$ of the superspace
${\cal M}/{\cal M}^2$, where $|\phi_i|=|f_i|=0$ for $1\leq i\leq
t$ and $|\phi_j|=|f_j|$ for $t+1\leq j\leq n$.

By induction on $l$ and using formula (3.1) of \cite{bk} we infer
that
$$\psi_1\ldots\psi_l(g_1\ldots g_l)=\sum_{1\leq i\leq l}(-1)^{\sum_{1\leq t < i}|g_i||g_t|+\sum_{1 < t}|g_i||\psi_t|}
\psi_1(g_i)\times$$
$$\sum_{\sigma(1)=i}(-1)^{\sum_{2\leq s< t, \sigma(s)
>\sigma(t)}|g_{\sigma(s)}||g_{\sigma(t)}|+\sum_{2\leq j < k}|\psi_k||g_{\sigma(j)}|}
\psi_2(g_{\sigma(2)})\ldots\psi_l(g_{\sigma(l)})$$ for
$\psi_1,\ldots ,\psi_l\in Lie(X)$ and $g_1,\ldots , g_l\in {\cal
M}$. This implies
$$\psi_1\ldots\psi_l(g_1\ldots g_l)=\sum_{\sigma\in S_l}
(-1)^{\sum_{s< t, \sigma(s)
>\sigma(t)}|g_{\sigma(s)}||g_{\sigma(t)}|+\sum_{j < i}|\psi_i||g_{\sigma(j)}|}
\psi_1(g_{\sigma(1)})\ldots\psi_l(g_{\sigma(l)}).$$ In particular,
$$\prod_{1\leq i\leq n}\phi_i^{s_i}(\prod_{1\leq i\leq
n}f_i^{d_i})=\pm\delta_{s_1 ,\ d_1}\ldots\delta_{s_n , \ d_n}s_1
!\ldots s_n ! ,$$ where $l=s_1 +\ldots +s_n=d_1+\ldots +d_n$ and
$s_{t+1},\ldots , s_n , d_{t+1},\ldots , d_n\in \{0, 1\}$.
Comparison of dimensions shows that $g_X$ is an isomorphism of
superalgebras.

By definition, $\delta_{Dist(X)}(\phi)(f_1\otimes f_2)=
\phi(f_1f_2)=\epsilon_X(f_1)f_2 +f_1\epsilon_X(f_2)$ for any
$\phi\in Lie(X)$, that is $\delta_{Dist(X)}(\phi)=\phi\otimes 1
+1\otimes\phi$. The same formula (3.1) of \cite{bk} implies that
$s_X(f)+f\in {\cal M}^2$ for $f\in {\cal M}$. In particular,
$s_{Dist(X)}(\phi)(f)= \phi(s_X(f))=-\phi(f)$ for $\phi\in
Lie(X)$, that is $s_{Dist(X)}(\phi)=-\phi$. Thus $g_X$ is an
isomorphism of filtered Hopf superalgebras.

For $A\in SAlg_K$ let $A[\varepsilon_0, \varepsilon_1]$ be a
(commutative) superalgebra of dual numbers. By definition,
$A[\varepsilon_0, \varepsilon_1]=\{a+\varepsilon_0 b
+\varepsilon_1 c| a, b, c\in A\}, |\varepsilon_i|=i,
\varepsilon_i\varepsilon_j=0, i, j\in \{0, 1\}$. We have two
morphism of superalgebras $p_A : A[\varepsilon_0,
\varepsilon_1]\to A$ and $i_A : A\to A[\varepsilon_0,
\varepsilon_1]$ defined by $a+\varepsilon_0 b +\varepsilon_1
c\mapsto a$ and $a\mapsto a$ respectively. Define the functor
${\bf Lie}(G)$ as
$${\bf Lie}(G)(A)=\ker (G(A[\varepsilon_0, \varepsilon_1])\stackrel{G(p_A)}{\to} G(A),
A\in SAlg_K .$$ It is called {\it Lie superalgebra functor} of
$G$. Let $V$ be a superspace. Define the functor $V_a$ from
$SAlg_K$ to the category of vector superspaces, by
$V_a(A)=V\otimes A$. The following lemma is obvious (see also
\cite{w, f1}).
\begin{lm}
There is an isomorphism of abelian group functors $Lie(G)_a\simeq
{\bf Lie}(G)$ given by
$$(v\otimes a)(f)=\epsilon_G(f)
+(-1)^{|a||f|}\varepsilon_{|v\otimes a|}v(f)a, v\in Lie(G)=({\cal
M}/{\cal M}^2)^*, a\in A, f\in K[G].$$
\end{lm}
If we identify $Lie(G)\otimes A$ with $Hom_K({\cal M}/{\cal M}^2 ,
A)$ via $(v\otimes a)(f)=(-1)^{|a||f|}v(f)a$, then the above
isomorphism can be represented as $$u\mapsto \epsilon_G
+\varepsilon_0 u_0 +\varepsilon_1 u_1, u\in Hom_K({\cal M}/{\cal
M}^2 , A).$$ Besides, this isomorphism induces the $A$-supermodule
structure on ${\bf Lie}(G)(A)$. The supergroup $G$ acts on the
functor ${\bf Lie}(G)$ by
$$(g, x)\mapsto G(i_A)(g)x G(i_A)(g)^{-1}, g\in G(A), x\in {\bf Lie}(G)(A),
A\in SAlg_K.$$ This action is called {\it adjoint} and denoted by
${\bf Ad}$.
\begin{lm}
The adjoint action is linear. In particular, it induces a
supergroup morphism $G\to GL(Lie(G))$.
\end{lm}
Proof. Let $u=v\otimes a\in Lie(G)\otimes A$ and $g\in G(A)$.
Denote the element $G(i_A)(g)$ by $\bar g$. Then
$${\bar g}u{\bar g}^{-1}(f)=\sum g(f_1)u(f_2)g(s_G(f_3))=
\epsilon_G(f)+\sum (-1)^{|a||f|+|f_1||v|}\varepsilon_{|u|}g(f_1)
v(f_2)g(s_G(f_3))a =$$
$$
({\bar g}(v\otimes 1){\bar g}^{-1}a)(f),
$$ where $(\delta_G\otimes
1)\delta_G(f)=(1\otimes\delta_G)\delta_G(f)=\sum f_1\otimes
f_2\otimes f_3, f\in K[G]$.
\begin{lm}
The differential of ${\bf Ad}$ coincides with $-{\bf ad}$.
\end{lm}
Proof. Denote $Lie(G)$ by $L$. We have a commutative diagram
$$\begin{array}{ccc}
G(K[\varepsilon_0, \varepsilon_1]) & \stackrel{Ad}{\to} &
GL(L)(K[\varepsilon_0, \varepsilon_1]) \\
\uparrow & & \uparrow \\
Lie(G) & \stackrel{d(Ad)}{\to} & gl(L).
\end{array}
$$
Notice that the image of $A\in gl(L)$ in $GL(L)(K[\varepsilon_0,
\varepsilon_1])$ is equal to $id_L +\varepsilon_0 A_0
+\varepsilon_1 A_1$. Choose $x, y\in Lie(G)$. Then
$$t=(d({\bf Ad})(x))(y)=(\epsilon_G +\varepsilon_{|x|} x)(\epsilon_G +\varepsilon'_{|y|} y)
(\epsilon_G -\varepsilon_{|x|} x),$$ where the product is computed
in $GL(L)((K[\varepsilon_0, \varepsilon_1])[\varepsilon'_0,
\varepsilon'_1])$. Further,
$$t(f)=\epsilon_G(f) + \varepsilon'_{|y|}(y(f) +\varepsilon_{|x|}(\sum (-1)^{|x||y|}
x(f_1)y(f_2)\epsilon_G(f_3)-\epsilon_G(f_1)y(f_2)x(f_3)))=$$
$$\epsilon_G(f) + \varepsilon'_{|y|}(y(f)-\varepsilon_{|x|}[y,
x](f)).$$

Following \cite{dg}, II, \S 4, we will denote the image of $u\in
{\bf Lie}(G)(A)$ in $G(A[\varepsilon_0, \varepsilon_1])$ by
$e^{\varepsilon_0 u_0 +\varepsilon_1 u_1}$.

\section{Quotient $K$-sheafs}

Let $G$ be a group dur $K$-sheaf and assume that $G$ acts freely
on a dur $K$-sheaf $X$, that is for any $R\in SAlg_K$ the group
$G(R)$ acts freely on $X(R)$. Then the functor $R\to
(X/G)_{(n)}(R)=X(R)/G(R)$ satisfies the properties (*) and (**)
for ff-coverings. The proof of this fact can be copied from
\cite{jan}, part I (5.5). Call the above functor $(X/G)_{(n)}$ a
{\it naive} quotient, the dur $K$-sheaf
$\tilde{\tilde{(X/G)_{(n)}}}$ a {\it quotient} dur $K$-sheaf (of
$X$ by $G$) and denote it by $\tilde{\tilde{X/G}}$. Then
$(X/G)_{(n)}\subseteq \tilde{\tilde{X/G}}$ and there is a
canonical $G$-invariant morphism $\tilde{\tilde{\pi}} : X\to
\tilde{\tilde{X/G}}$. Besides, for any other $G$-invariant
morphism of dur $K$-sheafs $h : X\to Z$ there is a unique morphism
$v : \tilde{\tilde{X/G}}\to Z$ such that $h=v\tilde{\tilde{\pi}}$
and these properties define $\tilde{\tilde{X/G}}$ uniquely up to
an isomorphism. Analogous statements are valid for $K$-sheafs and
there is an inclusion $\tilde{X/G}\subseteq \tilde{\tilde{X/G}}$.

Let $X$ be an affine superscheme and let $G$ be an affine
supergroup acting on $X$ via $f:X\times G \to X$. If $\tau :
K[X]\to K[X]\otimes K[G]$ is a comorphism dual to $f$, then
$K[X]^G=\{a\in K[X]|\tau(a)=a\otimes 1\}$ is a supersubalgebra of
$K[X]$. The embedding $K[X]^G\to K[X]$ induces a $G$-invariant
morphism  $i : X\to SSp \ K[X]^G$. In particular, there are
uniquely defined morphisms $i' : \tilde{X/G}\to SSp \ K[X]^G$ and
$i" : \tilde{\tilde{X/G}}\to SSp \ K[X]^G$ such that
$i'\tilde{\pi}=i=i"\tilde{\tilde{\pi}}$.
\begin{pr}
Suppose that $\tilde{\tilde{X/G}}$ (or $\tilde{X/G}$,
respectively) is an affine superscheme. Then $i"$ (or $i'$,
respectively) is an isomorphism and  $K[X]$ is a faithfully flat
$K[X]^G$-supermodule. If $\tilde{X/G}$ is an affine superscheme
and $G$ is algebraic, then $K[X]^G\leq K[X]$.
\end{pr}
Proof. Let $\tilde{\tilde{X/G}}=SSp \ R$ or $\tilde{X/G}=SSp \ R$.
There is a canonical isomorphism $pr_X\times f : X\times G\simeq
X\times_{SSp \ R} X$ (see \cite{jan}, part I (5.5)) that is dual
to the isomorphism of superalgebras $\phi : K[X]\otimes_R
K[X]\simeq K[X]\otimes K[G]$ defined by $a\otimes b\mapsto \sum
ab_1\otimes b_2$, where $\tau(b)=\sum b_1\otimes b_2$ and $a, b\in
K[X]$. Repeating the proof of \cite{jan}, part I (5.7) we obtain
an isomorphism of $B$-superalgebras $B\otimes_R K[X]\simeq
B\otimes K[G]$, where $B$ is either ff-covering of $R$ or $R\leq
B$. Using a superversion of Proposition 4 from \cite{bur}, I, \S 3
we see that $K[X]$ is a faithfully flat $R$-supermodule. Thus
$K[X]^G\subseteq R\subseteq K[X]$. If $G$ is algebraic, then
$B\otimes K[G]$ is a finitely presented $B$-superalgebra. A
superversion of Lemma 1.4, \cite{dg}, I implies that $R\leq K[X]$
in the case of $K$-sheafs. Composition of $\phi$ and an exact
sequence $$0\to R\to K[X]\begin{array}{c}\to \\
\to\end{array} K[X]\otimes_R K[X]$$ gives $R\subseteq K[X]^G$.
\begin{pr}
Assume that $R$ is a supersubalgebra of $K[X]^G$ and the canonical
morphism $X\times G\to X\times_{SSp \ R} X$ is an isomorphism. If
$K[X]$ is a faithfully flat $R$-supermodule, then $R=K[X]^G$ and
$\tilde{\tilde{X/G}}\simeq SSp \ R$. Additionally, if $R\leq
K[X]$, then $\tilde{X/G}\simeq SSp \ R$.
\end{pr}
Proof. Consider a $G$-invariant morphism of  dur $K$-sheafs $X\to
Z$. By Yoneda's lemma this morphism is uniquely defined by some
element $z\in Z(K[X])$. The $G$-invariance of this morphism is
equivalent to the following property. For any $A\in SAlg_K$ and
for arbitrary $\alpha\in X(A), \beta\in G(A)$ we have
$Z(\alpha)(z)=Z((\alpha\otimes\beta)\tau)(z)$. Set $A=K[X]\otimes
K[G]$ and $\alpha : a\mapsto a\otimes 1, \beta : b\mapsto 1\otimes
b$ for $a\in K[X]$ and $b\in K[G]$. Then
$\alpha\otimes\beta=id_{K[X]\otimes K[G]}$. Therefore $z$ belongs
to the kernel of
$$Z(K[X])\begin{array}{c}\stackrel{Z(\alpha)}{\to} \\
\stackrel{Z(\tau)}{\to}\end{array} Z(K[X]\otimes K[G]) .$$ On the
other hand, the above diagram can be identified with
$$Z(K[X])\begin{array}{c}\stackrel{Z(\gamma_1)}{\to} \\
\stackrel{Z(\gamma_2)}{\to}\end{array} Z(K[X]\otimes_R K[X])$$ via
bijection $Z(K[X]\otimes_R K[X])\to Z(K[X]\otimes K[G])$, where
$\gamma_1 : a\mapsto a\otimes 1$ and $\gamma_2 : a\mapsto 1\otimes
a$ for $a\in K[X]$. Setting $Z=A^{1|1}$ we obtain that
$K[X]^G\subseteq R$. Furthermore, if $Z$ is a dur $K$-sheaf, then
$z$ belongs to the image of the map $Z(R)\to Z(K[X])$ induced by
the inclusion $R\to K[X]$. The above quoted Yoneda's lemma
completes the proof. These arguments can be repeated per verbatim
for $K$-sheafs.

\section{Coideal supersubalgebras of Hopf superalgebras}

Let $A$ be a commutative Hopf superalgebra. Its supersubalgebra
$B\subseteq A$ is called a left (or right, respectively) {\it
coideal} iff $B$ is a left (or right, respectively)
$A$-supercomodule. A typical example is as follows. If $I$ is a
superideal and a coideal of $A$, then $C=A/I$ is a superbialgebra
and $B=A^C$ (or $B={^{C}}A$, respectively) is a left (or right,
respectively) coideal, see \cite{t2}. Denote by ${\cal H}^C$ the
category whose objects are right $A$-supermodules and
$C$-supercomodules simultaneously, together with even morphisms
such that $\tau_M(ma)=\sum (-1)^{|c_2||a_1|} m_1 a_1\otimes c_2
p(a_2)$, where $\delta_A(a)=\sum a_1\otimes a_2$ and
$\tau_M(m)=\sum m_1\otimes c_2$ for $m\in M, a\in A, M\in {\cal
H}^C$ and $p : A\to C=A/I$ is the canonical epimorphism. For
example, $A, C\in {\cal H}^C$.

Symmetrically, let $B$ be a left coideal supersubalgebra of $A$.
Denote by ${_{B}}{\cal H}$ the category whose objects are left
$B$-supermodules and $A$-supercomodules simultaneously, together
with even morphisms such that $\tau_M(bm)=\sum
(-1)^{|b_2||a_1|}b_1 a_1\otimes b_2 m_2$, where $\delta_A(b)=\sum
b_1\otimes b_2$ and $\tau_M(m)=\sum a_1\otimes m_2$ for $b\in B,
m\in M, M\in {_{B}}{\cal H}$. For example, $B, A\in {_{B}}{\cal
H}$.

\begin{lm}
Let $N$ be a right $A$-supermodule, $B$ be a coideal superalgebra
of $A$ and $M\in {_{B}}{\cal H}$. The linear map $\xi : N\otimes
M\to N\otimes M$, defined by $\xi(n\otimes m)=\sum na_1\otimes
m_2$, induces an isomorphism of superspaces $N\otimes_B M\to
N\otimes \overline{M}$, where $\overline{M}=M/B^{+}M$ and
$B^{+}=B\bigcap\ker\epsilon_A$.
\end{lm}
Proof. It is not difficult to see that $\xi$ is an isomorphism of
superspaces. The inverse $\xi^{-1}$ of $\xi$ is defined by
$\xi^{-1}(n\otimes m)= \sum n s_A(a_1)\otimes m_2$. We have
$$\xi(nb\otimes m-n\otimes bm)=\sum nba_1\otimes m_2 -\sum
(-1)^{|a_1||b_2|}nb_1 a_1\otimes b_2 m_2=$$
$$-\sum(-1)^{|a_1||b_2|}nb_1 a_1\otimes (b_2-\epsilon_A(b_2)) m_2\in
N\otimes B^+ M$$ and
$$\xi^{-1}(n\otimes bm)=\sum (-1)^{|a_1||b_2|}ns_A(b_1 a_1)\otimes
b_2 m_2 =-(\sum ns(a_1 b_1)b_2\otimes m_2-\sum ns(a_1 b_1)\otimes
b_2 m_2).$$
\begin{lm}
Let $p : A\to C$ be an epimorphism of superbialgebras. For any
$N\in {\cal H}^C, M\in A-scomod$ the previously defined map $\xi$
induces an isomorphism of superspaces $N^C\otimes M\simeq
N\square_C M$.
\end{lm}
Proof. It can be checked easily that $\xi(N^C\otimes M)\subseteq
N\square_C M$ and therefore it remains to show that
$\xi^{-1}(N\square_C M)\subseteq N^C\otimes M$. Let $\{m_i\}_{i\in
I}$ be a homogeneous basis of $M$ and $\tau_M(m_i)=\sum_{k\in I}
a_{ik}\otimes m_k$ for $i\in I$. Notice that $|a_{ij}|=|m_i|+|m_j|
\pmod 2$ for $i, j\in I$. The condition $\sum_{i\in I} n_i\otimes
m_i\in N\square_C M$ is equivalent to equalities
$\tau_N(n_i)=\sum_{k\in I} n_k\otimes p(a_{ki})$ for $i\in I$.
Applying $\xi^{-1}$ we see that all we have to check is
$u_k=\sum_{i\in I}n_is_A(a_{ik})\in N^C$ for all $k$. But this
follows from
$$\tau_N(u_k)=\sum_{i, t, l\in I}(-1)^{|a_{lt}||a_{lk}|}n_t s_A(a_{lk})\otimes
p(a_{ti}s_A(a_{il}))=$$$$\sum_{t, l\in
I}(-1)^{|a_{lt}||a_{lk}|}n_t s_A(a_{lk})\otimes
p(\epsilon(a_{tl}))=u_k\otimes p(1).$$

Observe that a superideal $AB^+$ is also a coideal. It follows
from $\delta_A(b)=b\otimes 1 +\sum a_1\otimes b_2$ for  $b\in B,
b_2\in B^+$. Define functors $\Phi : {_{B}}{\cal H}\to
\overline{A}-scomod$ and $\Psi : \overline{A}-scomod\to
{_{B}}{\cal H}$ by $\Phi(M)=\overline{M}, \ \Psi(N)=
A\square_{\overline{A}} N$. The functor $\Psi$ is right adjoint of
$\Phi$ by \cite{t2}. Adjunctions $u_M : M\to
\Psi\Phi(M)=A\square_{\overline{A}}\overline{M}, \ v_N :
\Phi\Psi(N)=\overline{A\square_{\overline{A}} N}\to N$ are defined
by $m\mapsto\sum a_1\otimes \overline{m_2}$ and $\overline{\sum
n\otimes a}\mapsto\sum\epsilon_A(a)n$. Symmetrically, one can
define functors $\Theta : smod-B\to {\cal H}^C$ and $\Omega :
{\cal H}^C\to smod-B$ by $\Theta(M)=M\otimes_B A, \
\Omega(N)=N^C$. The functor $\Theta$ is left adjoint of $\Omega$
by \cite{t2}. The corresponding adjunctions $f_M : M\to
\Omega\Theta(M)=(M\otimes_B A)^C$ and $g_N :
\Theta\Omega(N)=N^C\otimes_B A\to N$ are defined by
$f_M(m)=m\otimes 1$ and $g_N(n\otimes a)=na$.
\begin{lm}
Let $V\in smod-A$ be a flat $B$-supermodule. Then $V\otimes_B
\Psi(N)$ is isomorphic to $V\otimes N$ via $v\otimes (\sum
a\otimes n)\mapsto\sum va\otimes n$, where $v\in V, \sum a\otimes
n\in \Psi(N)$.
\end{lm}
Proof. Maps $(id_A\otimes ?)\delta_A\otimes id_N$ and
$id_A\otimes\tau_N$ are morphisms of left $B$-supermodules. They
combine to a canonical isomorphism $V\otimes_B \Psi(N)\simeq
(V\otimes_B A)\square_{\overline{A}} N$ (see also Proposition 1.3,
\cite{t3}). On the other hand, $\xi : V\otimes_B A\to
V\otimes\overline{A}$ is an isomorphism of right
$\overline{A}$-supercomodules and its composition with $\xi\otimes
id_N$ gives $V\otimes_B \Psi(N)\simeq
V\otimes\overline{A}\square_{\overline{A}} N\simeq V\otimes N$.
\begin{lm}
Let $V\in A-scomod$ and $V$ be a coflat $C$-supercomodule. For any
$B$-supermodule $T$ there is an isomorphism $T\otimes V\simeq
\Theta(T)\square_C V$ defined by $t\otimes v\mapsto \sum t\otimes
c_1\otimes v_2$, where $\tau_V(v)=\sum c_1\otimes v_2$.
\end{lm}
Proof. Denote by $p : T\otimes A\to T\otimes_B A$ the canonical
epimorphism of superspaces. A sequence $0\to\ker p\to T\otimes
A\stackrel{p}{\to} T\otimes_B A$ is exact in the category
$C-scomod$. In particular, it induces a canonical isomorphism
$T\otimes_B (A\square_C V)\to (T\otimes_B A)\square_C V$.
Composition with $id_T\otimes\xi$ gives the isomorphism of the
lemma.
\begin{lm}
Let $A$ be a superalgebra and let $\phi : M\to N$ be an
epimorphism of free $A$-supermodules of the same finite rank. Then
$\phi$ is an isomorphism.
\end{lm}
Proof. Denote the rank of both $M$ and $N$ by $(r, s)$. Choose
free generators $m_i, n_i$ for $1\leq i\leq r+s=n$ of supermodules
$M$ and $N$ correspondingly. Additionally, assume that
$|m_i|=|n_i|=0$ if $1\leq i\leq r$ and $|m_i|=|n_i|=1$ otherwise.
Denote $\phi(m_i)=n'_i$ for $1\leq i\leq n$ and write
$n'_i=\sum_{1\leq k\leq n}b_{ik}n_k, n_i=\sum_{1\leq k\leq
n}c_{ik} n'_k$ for $1\leq i\leq n$. Matrices $B=(b_{ik})$ and
$C=(c_{ik})$ belong to $Mat_{r, s}(A)_0$ and $CB=E_n$. Considering
the last matrix equation modulo the superideal $AA_1$, we infer
that even diagonal blocks of the matrix $C$ are invertible. By
Lemma 1.7.2 of \cite{l}, (see also Theorem 3.1, \cite{ber}) $C$ is
invertible. In particular, the elements $n'_i$ are free generators
over $A$.
\begin{pr}\label{5.1}
If $M\in {\cal H}_B$, then a left $B$-module $M_{B^+}$ is flat.
\end{pr}
Proof. We follow the same ideas as in the proofs of Lemma 3.1 and
of Corollary 3.3 in \cite{mw}. First of all, it is enough to
consider the case when $M$ is finitely generated. Using the
isomorphism $K[G]\otimes_B M\to K[G]\otimes \overline{M}$ we see
that a left $K[G]$-supermodule $K[G]\otimes_B M$ is free of rank
$(m, n)$, where $m=\dim\overline{M}_0$ and $n=\dim\overline{M}_1$.
It is clear that there is an isomorphism of $B$-supermodules $f :
B^m\bigoplus (B^c)^n\to M$ which induces an isomorphism modulo
$B^+$. On the other hand, $B_{B^+}$ is a local algebra.  Thus $rad
M_{B^+}=(B^+ B_{B^+}) M_{B^+}$ and the localization of $f$,
denoted by $f_{B^+} : B_{B^+}^m\bigoplus (B_{B^+}^c)^n\to
M_{B^+}$, is an epimorphism (cf. \cite{kash}, 9.3.5, 9.4.1).
Tensoring with $K[G]_{B^+}$ produces an epimorphism $g :
K[G]_{B^+}^m\bigoplus (K[G]_{B^+}^c)^n\to K[G]_{B^+}\otimes
_{B_{B^+}}M_{B^+}$. The $K[G]_{B^+}$-supermodule
$K[G]_{B^+}\otimes _{B_{B^+}}M_{B^+}$ can be naturally identified
with $K[G]_{B^+}\otimes_B M\simeq K[G]_{B^+}\otimes\overline{M}$.
By Lemma 5.5, the map $g$ is an isomorphism and its composition
with $B_{B^+}^m\bigoplus (B_{B^+}^c)^n\to (K[G]_{B^+}^c)^n\to
K[G]_{B^+}$ is identified with $f_{B^+}$.
\begin{lm}
A Hopf superalgebra $H$ is a direct union of all of its finitely
generated Hopf supersubalgebras.
\end{lm}
Proof. Any finite subset $X\subseteq H$ is contained in a finite
dimensional supersubcoalgebra $C$. Let $c_i$ for  $1\leq i\leq t$
be a homogeneous basis of $C$. The supersubalgebra generated by
$c_i$ and $s_H(c_i)$ for $1\leq i\leq t$ is obviously a Hopf
superalgebra containing $X$.
\begin{ter}
The Hopf superalgebra $A$ is a flat left $B$-module.
\end{ter}
Proof. Without a loss of generality one can assume that $K$ is
algebraically closed and by Lemma 5.6 one can suppose that $A$ is
finitely generated. For any ${\cal M}\in Max(A)$, the map
$T_{{\cal M}} : A\to A$ defined by $T_{{\cal M}}(f)=\sum
\pi_{{\cal M}}(f_1)f_2$, where $\delta_A(f)=\sum f_1\otimes f_2$
and $\pi_{{\cal M}} : A\to A/{\cal M}=K$, is an automorphism of
the superalgebra $A$ because $A_1\subseteq {\cal M}$. The inverse
of $T_{{\cal M}}$ is $T_{s_A ({\cal M})}$ because $\pi_{s_A({\cal
M})}(f)=\pi_{{\cal M}}(s_A(f))$ for $f\in A$.

Next, the map $T_{{\cal M}}$ takes $B$ to $B$ and ${\cal M}$ to
${\cal M}_1=\ker\epsilon_A$. Using Proposition 5.1, Lemma 1.2 and
the left-hand side version of Proposition 8, I, \S 2 of \cite{bur}
we see that $A_{B^+}\otimes_{A} A_{{\cal M}_1}$ is a flat left
$B$-module. The canonical epimorphism of left $B$-modules
$A_{B^+}\otimes_{A} A_{{\cal M}_1}\to A_{{\cal M}_1}$ is split.
Therefore $A_{{\cal M}_1}$ is a flat left $B$-module. The
isomorphism $T_{{\cal M}}$ takes this module to $A_{{\cal
M}}^{\alpha_{{\cal M}}}$, where $\alpha_{{\cal M}}=T_{{\cal M}}
|_{B}$. We complete the proof by combining Proposition 1.1, the
remark after Corollary 1.2 and the same reductions as in \cite{mw}
(see also Lemma 7.1, III, \S 3 of \cite{dg}).
\begin{cor}
If $B$ is a Hopf supersubalgebra of $A$, then $A$ is a faithfully
flat (left and right) $B$-module.
\end{cor}
Proof. By Theorem \ref{5.1} it remains to show that ${\cal M}A\neq
A$ for all ${\cal M}\in Max(B)$. Since $J=A A_1$ is a Hopf
superideal and $J\bigcap B\subseteq {\cal M}$, it suffices to
check that ${\cal M}(A/J)\neq A/J$. On the other hand, $A/J$ is a
Hopf algebra and $B/J$ is its Hopf subalgebra. By Takeuchi's
theorem (see \cite{t2}), $A/J$ is a faithfully flat $B/J$-module.
The right-hand side statement can be deduced from the right-hand
side version of Theorem \ref{5.1}.

Let $G$ be an affine supergroup and let $H$ be its closed
supersubgroup. The supergroup $H$ is called {\it faithfully exact}
in $G$ if $K[G]$ is a faithfully exact right $K[H]$-supercomodule.
We identify the category of right $H$-supermodules $mod-H$ with
$K[H]-scomod$. In particular, the functor $\Psi=K[G]\square_{K[H]}
?$ can be identified with the inducing functor $ind^G_H$ (see
\cite{z}).
\begin{pr}
The quotient $\tilde{\tilde{G/H}}$ is affine iff $I_H=K[G]R^+$ and
$K[G]$ is a faithfully flat $R$-supermodule, where $R$ is a left
coideal supersubalgebra of $K[G]^H$. If $G$ is algebraic, then
$\tilde{G/H}$ is affine iff $I_H=K[G]R^+$ and $R\leq K[G]$.
\end{pr}
Proof. Lemma 5.1 and the equality $I_H=K[G]R^+$ guarantees that
$G\times H\simeq G\times_{SSp \ R} G$. Conversely, assume that the
quotient $\tilde{\tilde{G/H}}$ (or $\tilde{G/H}$) is affine. In
the following commutative diagram of superalgebras
$$\begin{array}{ccc}
K[G]\otimes_R K[G] & \stackrel{\xi}{\longrightarrow} & K[G]\otimes\overline{K[G]} \\
\searrow & & \swarrow \\
& K[G]\otimes K[H]  &
\end{array}$$
the horizontal and the left diagonal arrows are isomorphisms. Thus
$\overline{K[G]}\to K[H]$ is an isomorphism and it remains to
refer to Proposition 4.2.

\begin{ter}
The following statements are equivalent : \\
i) Quotient $\tilde{\tilde{G/H}}$ is affine; \\ ii) $H$ is a
faithfully exact supersubgroup of $G$; \\ iii) $ind^G_H$ induces
an equivalence of $mod-H$ with the full subcategory ${_{R}}{\cal
H}$;
\\ iv) $K[G]$ is an injective cogenerator in the category $mod-H$.
\end{ter}
Proof.  Combining Lemmas 5.1 -- 5.4 and Proposition 5.2 with the
proof of Theorems 1 and 2 from\cite{t2} we easily obtain that the
properties i), ii) and iii) are equivalent to each other. The
equivalence of iii) and iv) was proved in Lemma 1.7.

\section{Quotients by normal supersubgroups}

A group $K$-subfunctor $H$ of a group $K$-functor $G$ is said {\it
normal} if $H(A)\unlhd G(A)$ for all $A\in SAlg_K$. If $G$ is a
dur $K$-sheaf (or $K$-sheaf, respectively) and $H$ is its normal
subsheaf, then Lemma 2.3 implies that $\tilde{\tilde{G/H}}$ (or
$\tilde{G/H}$, respectively) is again group dur $K$-sheaf (or
group $K$-sheaf, respectively) and the canonical morphism
$G\to\tilde{\tilde{G/H}}$ (or $G\to\tilde{G/H}$, respectively) is
a morphism of group $K$-functors.

If $G$ is an affine supergroup and $H$ is its closed
supersubgroup, then $H\unlhd G$ iff it satisfies one of the
following conditions
$$\nu_r(f)=\sum (-1)^{|f_1||f_2|}f_2\otimes f_1 s_G(f_3)\in I_H\otimes K[G],$$
or
$$\nu_l(f)=\sum (-1)^{|f_1||f_2|}f_2\otimes s_G(f_1)f_3\in I_H\otimes K[G],$$
for any $f\in I_H$. The first condition is called {\it right
normality} condition and the second one is called {\it left
normality} condition. These conditions are different in general
(say for quantum groups) but for supergroups they are equivalent
because $s_G$ is an automorphism of superalgebra $K[G]$ of order
two (see \cite{pw}, 1.5, and \cite{z}, \S 2). The morphism $\nu_l
$ is a superalgebra morphism, it is dual to the morphism of
superschemes $G\times G\to G$ given by $(g_1, g_2)\mapsto
g_2^{-1}g_1 g_2$ for $g_1, g_2\in G(A)$ and $A\in SAlg_K$ (which
defines a right action of $G$ on itself by conjugations).
Symmetrically, $\nu_r$ is dual to the morphism $(g_1, g_2)\mapsto
g_2 g_1 g_2^{-1}$. From now on, all group $K$-functors are assumed
affine and all group subfunctors are assumed closed unless
otherwise stated.
\begin{lm}
If $H\unlhd G$, then $Lie(H)$ is a Lie superideal of $Lie(G)$.
\end{lm}
Proof. It is clear that ${\bf Ad}(G)\subseteq
Stab_{GL(Lie(G))}(Lie(H))$. Lemma 3.4 concludes the proof.

Let $H$ be a supersubgroup of $G$. Denote the normalizer of $H$ in
$G$ by $N_G(H)$. By definition, $$N_G(H)(A)=\{g\in G(A)|\mbox{for
any} \ A-\mbox{superalgebra} \ B \
G(\iota^B_A)(g)H(B)G(\iota^B_A)(g)^{-1}=H(B) \}.$$
\begin{pr}
The normalizer $N_G(H)$ is a closed supersubgroup.
\end{pr}
Proof. Consider $f\in I_H$. Let $\nu_l(f)=\sum u_1\otimes u_2$
modulo $I_H\otimes K[G]$, where $u_1$ are linearly independent
modulo $I_H$. The elements $u_2$ are called {\it companions} of
$f$. Let $I$ be a superideal of $K[G]$ generated by companions of
all $f$. Set $B=K[H]\otimes A$ and $\iota^B_A=1\otimes id_A$. Then
$h=p\otimes 1\in H(B)$, where $p$ is the canonical epimorphism
$K[G]\to K[H]$. We have
$$(G(\iota^B_A)(g)^{-1}hG(\iota^B_A)(g))(f)=\sum u_1\otimes
g(u_2)=0 .
$$
Thus $g(u_2)=0$. Conversely, if all $g(u_2)=0$, then
$G(\iota^B_A)(g)H(B)G(\iota^B_A)(g)^{-1}\subseteq H(B)$ for any
$A$-superalgebra $B$. It implies $N_G(H)=V(I+s_G(I))$.
\begin{lm}
If $H\unlhd G$, then $R={^{H}}K[G]=K[G]^H$ is a Hopf
supersubalgebra of $K[G]$.
\end{lm}
Proof. It is enough to observe that $f\in R=K[G]^H$ iff
$f(gh)=f(g)$ for all $g\in G(A), h\in H(A)$ and $A\in SAlg_K$. In
particular, $f(hg)=f(gg^{-1}hg)=f(g)$ implies  that the antipode
$s_G$ induces an automorphism of $R$.

If $\phi : G\to H$ and $\psi : L\to H$ are morphisms of affine
supergroups, then $G\times_H L$ is a supersubgroup  of  $G\times
L$. In the case $L\leq H$ we can identify the fiber product
$G\times_H L$ with a supersubgroup of $G$ that we call a {\it
preimage} of $\phi^{-1}(L)$. Besides,
$K[\phi^{-1}(L)]=K[G]/K[G]\phi^*(I_L)$. In particular, a kernel
$N=\ker\phi\leq G$ coincides with $G\times_H E$, where $E$ is the
trivial supersubgroup of $H$. It is also a kernel in the category
of supergroups (see Proposition 1.6.1 from \cite{pw}). As in
\cite{pw}, we have that $N \unlhd G$ and, moreover, $N\unlhd
\phi^{-1}(L)$.
\begin{ter}
A supergroup $\tilde{\tilde{G/N}}$ is isomorphic to a
supersubgroup of $H$ (the image of $\phi$). More precisely,
$Im\phi= SSp \ K[H]/\ker\phi^* =SSp \ Im\phi^*$. Thus $N$ is
faithfully exact supersubgroup of $G$. If $G$ and $H$ are
algebraic, then $\tilde{G/N}=\tilde{\tilde{G/N}}$.
\end{ter}
Proof. The definition of $I_N$  guarantees that the canonical
morphism $G\times N\to G\times_{Im\phi} G$ is an isomorphism.
Besides, $Im\phi^*\subseteq K[G]^N$. Proposition 5.2 and Corollary
5.1 imply the first statement. The second statement is deduced
easily from Lemma 1.1.
\begin{cor}
The canonical epimorphism $\phi^{-1}(L)\to L\bigcap Im\phi$, which
is dual to the embedding $K[H]/(I_L +\ker\phi^*)\to
K[G]/K[G]\phi^*(I_L)$, induces an isomorphism
$\tilde{\tilde{\phi^{-1}(L)/N}}\simeq L\bigcap Im\phi$. Moreover,
the diagram
$$\begin{array}{ccc}
\tilde{\tilde{G/N}} & \simeq & Im\phi \\
\uparrow & & \uparrow \\
\tilde{\tilde{\phi^{-1}(L)/N}} & \simeq & L\bigcap Im\phi
\end{array} ,$$
where the vertical maps are natural embeddings, is commutative.
\end{cor}
\begin{pr}
The following statements are equivalent : \\
i) A quotient $\tilde{\tilde{G/H}}$ is affine for any algebraic
supergroup $G$
and for any normal supersubgroup $H$ of $G$; \\
ii) Quotient $\tilde{\tilde{G/H}}$ is affine for any affine
supergroup $G$ and for any normal supersubgroup $H$ of $G$.
\end{pr}
Proof. We have to check only the implication i) $\rightarrow$ ii).
According to Lemma 5.6, $K[G]$ is a direct union of its finitely
generated Hopf supersubalgebras, say $K[G]=\bigcup_{\alpha\in
{\cal A}}B_{\alpha}$. Set $I_{\alpha}=B_{\alpha}\bigcap I_H$ for
$\alpha\in {\cal A}$. By Theorem 5.1 and by i) for any pair
$H_{\alpha}=SSp \ B_{\alpha}/I_{\alpha}\unlhd G_{\alpha}=SSp \
B_{\alpha}$ we obtain that $B_{\alpha}=K[G_{\alpha}]$ is a
faithfully flat (left and right)
$R_{\alpha}=K[G_{\alpha}]^{H_{\alpha}}$-module and
$I_{\alpha}=B_{\alpha}R_{\alpha}^+$. It is clear that
$\bigcup_{\alpha\in {\cal A}}R_{\alpha}=R=K[G]^H$. By Lemma 7.1,
III, \S 3 of \cite{dg}, $K[G]$ is faithfully flat (left and right)
$R$-module. It remains to observe that $I_H=\bigcup_{\alpha\in
{\cal A}}I_{\alpha}=\bigcup_{\alpha\in {\cal
A}}B_{\alpha}R_{\alpha}^+ = K[G]R^+$.
\begin{rem}
If $G$ is algebraic, $H\unlhd G$ and $\tilde{\tilde{G/H}}$ is
affine, then $\tilde{\tilde{G/H}}=\tilde{G/H}$. In fact,
$I_H=K[G]R^+$, where $R=K[G]^H$. By Lemma 1.1 the superideal $I_H$
is finitely generated. Moreover, it is generated by some finite
subset from $R^+$. By Lemma 5.6 this subset is contained in a
finitely generated Hopf supersubalgebra $B\subseteq R$. Using
Propositions 5.2 and 4.2 we see that $B=R$.
\end{rem}

\begin{pr}
Let $G$ be an algebraic supergroup and $H\leq G$. There is a
linear representation $\phi : G\to GL(V)$ such that $\phi :
G\simeq Im \phi$ and $\phi|_H : H\simeq Stab_G(W)$ for a suitable
supersubspace $W\subseteq V$.
\end{pr}
Proof. There is a finite dimensional supersubcomodule $V\subseteq
K[G]$ containing all generators of $K[G]$ as well as all
generators of $I_H$. Let $v_1,\ldots , v_{m+n}$ be a basis of $V$
such that $|v_i|=0$ for $1\leq i\leq m$ and $|v_i|=1$ otherwise.
Additionally, assume that $v_1,\ldots, v_s, v_{m+1},\ldots ,
v_{m+t}$ for $s\leq m$ and $t\leq n$ is a basis of $W=I_H\bigcap
V$.  We have a morphism of supergroup $\phi : G\to GL(V)$ defined
by $g\mapsto (g(r_{ij}))$ for $g\in G(A)$ and $A\in SAlg_K$, where
$\tau_V(v_i)=\delta_G(v_i)=\sum_{1\leq j\leq m+n} v_j\otimes
r_{ji}$ for $1\leq i\leq m+n$. If $M=\{1,\ldots, s, m+1,\ldots
m+t\}$ and $i\in M, j\not\in M$, then $r_{ji}\in I_H$.
Superalgebra $Im\phi^*$ is generated by the elements $r_{ij}$ and
by the multiplicative set generated by determinants of even blocks
of the matrix $(r_{ij})$. On the other hand, $v_i=\sum_{1\leq
j\leq m+n}\epsilon_G(v_j)r_{ji}$ for every $i$. Thus $Im\phi^*
=K[G]$ and $\ker\phi=E$. Finally, if $i\in M$, then
$v_i=\sum_{j\not\in M}\epsilon_G(v_j)r_{ji}$. In other words,
$\phi(H)=Stab_{G}(W)$.

\begin{pr}
Let $G$ be a group $K$-sheaf and $N_1 \leq N_2$ be group subsheafs
of $G$. If $N_1\unlhd G$, then $H=\tilde{N_2/N_1}$ is a group
subsheaf of $M=\tilde{G/N_1}$ and
$\tilde{M/H}\simeq\tilde{G/N_2}$. Additionally, if $N_2\unlhd G$,
then $H\unlhd M$ and the last isomorphism is an isomorphism of
group sheafs. Analogous statements are valid for dur $K$-sheafs.
\end{pr}
Proof. It is an easy consequence of the universal property of
quotients combined with Lemma 2.3.

Now we can formulate and prove the main result of this paper.
\begin{ter}
If $G$ is an affine supergroup and $N$ is a normal supersubgroup
of $G$, then $\tilde{\tilde{G/N}}$ is again an affine supergroup.
\end{ter}

Let $G$ be an affine supergroup and $N$ is a normal supersubgroup
of $G$. By Proposition 6.2  one can assume that $G$ is algebraic.
Define the supersubgroup $\overline{N}\leq G$ in such way that
$I_{\overline{N}}=K[G]R^+$, where $R=K[G]^N$. By Lemma 6.2 we have
$N\leq \overline{N}\unlhd G$ and
$\tilde{\tilde{G/\overline{N}}}\simeq SSp \ R$. It remains to
prove that $N=\overline{N}$.
\begin{lm}
The superalgebra $K[\overline{N}]^N$ coincides with $K$.
\end{lm}
Proof. The canonical isomorphism $K[G]\otimes_R K[G]\to
K[G]\otimes K[\overline{N}]$ is an isomorphism of right
$K[N]$-supercomodules. Consider the following exact sequence of
$R$-supermodules
$$0\to R\to K[G]\stackrel{\phi}{\to} K[G]\otimes K[N] ,$$
where $\phi(a)=\sum a_1\otimes\overline{a_2}-a\otimes\overline{1}$
for $a\in K[G]$ and $\delta_G(a)=\sum a_1\otimes a_2$. Using
Theorem 5.1 we obtain $K[G]\otimes_R R =(K[G]\otimes_R
K[G])^N\simeq K[G]\otimes K[\overline{N}]^N$ which implies
$K[\overline{N}]^N=K$.

From now on we assume that $G$ is algebraic and $K[G]^N=K$ unless
stated otherwise. Without a loss of generality one can assume that
$K$ is algebraically closed. Up to the end of this section $char
K=p > 0$.

The radical $r$ of the superalgebra $K[G]$ is a Hopf superideal.
In fact, a superalgebra $K[G]/r\otimes K[G]/r$ is reduced as the
coordinate algebra of an affine variety $Max(K[G])\times
Max(K[G])$. A supergroup $G_{red}$ corresponding to the Hopf
superideal $r$ is pure even. In other words, it is an affine group
(= affine group scheme). Besides, $G_{red}\leq G_{ev}$, where
$I_{G_{ev}}=K[G]K[G]_1$. It is clear that $U=N\bigcap
G_{red}\unlhd G_{red}$.
\begin{lm}
We have $K[G_{red}]^{N_{red}}=K$. In particular, $G_{red}\leq U$
or equivalently, $U=N_{red}$.
\end{lm}
Proof. Assume that an element $f\in K[G]_0$ represents a
$N_{red}$-invariant in $K[G_{red}]$, that is $\delta_G(f)-f\otimes
1\in r\otimes K[G] + K[G]\otimes (r+I_N)$. Since $r$ is a
nilpotent ideal, it follows that for a sufficiently large integer
$M > 0$ we have $\delta(f^{p^M})-f^{p^M}\otimes 1\in K[G]\otimes
I_N$. In particular, $f^{p^M}=a\in K$. Thus $f=b+x$, where $b\in
K, b^{p^M}=a$ and $x\in r$. The second statement follows from
Theorem 4.3 of \cite{t1}.

According to Proposition 6.3 we can write $G\leq GL(V)$ and
$N=Stab_G(W)$ for suitable supersubspaces $W\subseteq V$. Using
the notations from Proposition 6.3 one can depict a matrix from
$GL(V)(F)$ for $F\in SAlg_K$ as
$$\left(\begin{array}{cc}
A & B \\
C & D \end{array}\right) ,$$ where $A=(A_{ij})\in GL_m(F_0) , \
B=(B_{ij})\in M_{m\times n}(F_1), \ C=(C_{ij})\in M_{n\times
m}(F_1)$ and $D=(D_{ij})\in GL_n(F_0)$ for $i, j=1, 2,$ where the
blocks $A_{ij}$ and $B_{ij}$ ($C_{ij}$ and $D_{ij}$, respectively)
have $s$ rows if $i=1$, and $m-s$ rows if $i=2$ ($t$ rows if
$i=1$, and $n-t$ rows if $i=2$, respectively). Symmetrically, the
blocks $A_{ij}$ and $C_{ij}$ ($B_{ij}$ and $D_{ij}$, respectively)
have $s$ columns if $j=1$, and $m-s$ columns if $j=2$ ($t$ columns
if $j=1$, and $n-t$ columns if $j=2$, respectively). The
supergroup $S(F)= (Stab_{GL(V)}(W))(F)$ consists of all matrices
with $A_{21}=0, B_{21}=0, C_{21}=0$ and $D_{21}=0$.

The open subfunctor $GL(V)_{f}$, where
$f=\det(A_{11})\det(D_{11})$, contains $S$. Let $U$ be a closed
supersubscheme of $GL(V)_{f}$ defined by equations $A_{11}=E_{s},
A_{22}=E_{m-s}, D_{11}=E_t, D_{22}=E_{n-t}, B_{11}=0, B_{12}=0,
B_{22}=0, C_{11}=0, C_{12}=0, C_{22}=0, D_{12}=0$ and $A_{12}=0$.
It is obvious that $U$ is an (unipotent) supersubgroup of $GL(V)$.
\begin{pr}
We have an isomorphism of superschemes $\psi : GL(V)_f\simeq
U\times S$ commuting with the right action of $S$ given by
multiplication.
\end{pr}
Proof. The above isomorphism $\psi : GL(V)_f\to U\times S$ is
defined by the rule
$$\left(\begin{array}{cc}
A & B \\
C & D \end{array}\right)\mapsto \left(\begin{array}{cc}
A' & B' \\
C' & D' \end{array}\right)\times \left(\begin{array}{cc}
A" & B" \\
C" & D" \end{array}\right),$$ where $$\left(\begin{array}{cc}
A' & B' \\
C' & D' \end{array}\right)\in U, \left(\begin{array}{cc}
A" & B" \\
C" & D" \end{array}\right)\in S,$$ and
$$A_{11}=A"_{11}, \ A_{12}=A"_{12}, \ B_{11}=B"_{11}, \
B_{12}=B"_{12},$$ $$C_{11}=C"_{11}, \ C_{12}=C"_{12}, \
D_{11}=D"_{11}, \ D_{12}=D"_{12},$$
$$\left(\begin{array}{cc}
A'_{21} & B'_{21} \end{array}\right)=\left(\begin{array}{cc}
A_{21} & B_{21}
\end{array}\right)\left(\begin{array}{cc}
A_{11} & B_{11} \\
C_{11} & D_{11} \end{array}\right)^{-1},$$
$$\left(\begin{array}{cc} C'_{21} & D'_{21}
\end{array}\right)=\left(\begin{array}{cc} C_{21} & D_{21}
\end{array}\right)\left(\begin{array}{cc}
A_{11} & B_{11} \\
C_{11} & D_{11} \end{array}\right)^{-1} ,$$
$$A"_{22}=A_{22}-A'_{21}A_{12}-B'_{21}C_{12},
B"_{22}=B_{22}-A'_{21}B_{12}-B'_{21}D_{12},$$
$$C"_{22}=C_{22}-C'_{21}A_{12}-D'_{12}C_{12},
D"_{22}=D_{22}-C'_{21}B_{12}-D'_{12}D_{12}.$$ The inverse morphism
is just the multiplication map.
\begin{lm}
The supergroup $G$ is a closed supersubscheme of $GL(V)_f$.
\end{lm}
Proof. Since $N\subseteq GL(V)_f$, we see that the image of $f$ in
$K[G]$ is invertible modulo the superideal $I_N$. On the other
hand, $I_N\subseteq r$ and therefore, $f\in K[G]^*$.
\begin{pr}
The supergroup $G$ coincides with $N$.
\end{pr}
Proof. The naive quotient morphism $GL(V)_f\to GL(V)_f/S$ can be
identified with the composition of $\psi$ and the projection
$U\times S\to U$. In particular,
$\tilde{\tilde{GL(V)_f/S}}=GL(V)_f/S\simeq U$. The induced
morphism $\pi : G\to GL(V)_f/S$ is dual to the composition of the
embedding $K[U]\otimes 1\to K[U]\otimes K[S]$ and the epimorphism
$K[U]\otimes K[S]\to K[G]$. The last epimorphism is a morphism of
right $K[N]$-supercomodules. Since $K[U]\otimes 1 = (K[U]\otimes
K[S])^S$ it follows that $Im\pi^*= K$, that is $G\subseteq S$.

\section{Quotients of finite supergroups}

Let $G$ be an affine supergroup and let $H, N$ be supersubgroups
of $G$ such that $H$ normalizes $N$. Denote the semi-direct
product of $H$ and $N$ by $H\dot{\times} N$. More precisely,
$(H\dot{\times} N)(A)=H(A)\times N(A)$ and $(h, n)(h', n')= (hh' ,
h'^{-1}nh'n')$ for any $h, h'\in H(A), n, n'\in N(A), A\in
SAlg_K$. We have a natural morphism $ g: H\dot{\times} N\to G, (h,
n)\mapsto hn$.  By Theorem 6.1 the image of $g$ is a closed
supersubgroup of $G$ which is denoted by $HN$. The construction in
Section 2 yields for any $A\in SAlg_K$
$$HN(A)=\{g\in G(A)| \mbox{there exists} \ \mbox{ff-covering} \ B \ \mbox{of} \ A \ \mbox{such that}
\ G(\iota^B_A)(g)\in H(B)N(B)\}.$$ For example, assume that $N$ is
a kernel of an epimorphism $\pi : G\to L$. The preimage
$\pi|_H^{-1}(\pi|_H(H))$ coincides with $HN$. Theorem 6.1 and
Corollary 6.1 imply that $I_{HN}=K[G](\pi^*(\ker\epsilon_L)\bigcap
I_H)$.
\begin{pr}(see \cite{jan}, part I, (6.2)) The quotient $\tilde{\tilde{HN/H}}$
is isomorphic to  $\tilde{\tilde{N/(N\bigcap H)}}$.
\end{pr}
Proof. The image of the canonical inclusion $(N/(N\bigcap
H))_{(n)}\to (HN/H)_{(n)}$ is dense with respect to the
Grothendieck topology of  ff-coverings. In fact, if $gH(A)\in
HN(A)/H(A)$, then there is a ff-covering $B$ of $A$ such that
$g'H(B)=h'' n''H(B)=h'' n'' (h'')^{-1}H(B)$, where
$g'=G(\iota^B_A)(g), h''\in H(B), n''\in N(B)$. Therefore,
$$\tilde{\tilde{N/(N\bigcap H)}}=\tilde{\tilde{(N/(N\bigcap
H))_{(n)}}}=\tilde{\tilde{(HN/H)_{(n)}}}=\tilde{\tilde{HN/H}}.$$

\begin{lm}
Let $A$ be a  finitely generated (commutative) superalgebra and
$I$ be a nilpotent superideal of $A$. If $\dim A/I$ is finite,
then $A$ is finite dimensional.
\end{lm}
Proof. Since $A_1$ is a finitely generated $A_0$-module, all we
have to check is that $\dim A_0 <\infty$. Denote by $V$ a finite
dimensional subspace of $A_0$ such that $V+I_0=A_0$. Choose a
non-negative integer $k$ such that $I^{k+1}=0$. Let
$I_0=\sum_{1\leq i\leq l}A_0z_i$. For any $a\in A_0$ we have
$$a=v+\sum_{1\leq i\leq l}a_i z_i, \mbox{ where } v\in V \mbox{ and } a_i\in A_0 .$$
Repeating this procedure for the coefficients $a_i$ we obtain
$$a=\sum_{0\leq t\leq lk}\sum_{1\leq i_1, \ldots,  i_t\leq
l}v_{i_1 ,\ldots , i_t}z_{i_1}\ldots z_{i_t},$$ where $v, v_{i_1
,\ldots , i_t}\in V.$ Therefore $\dim
A_0\leq\frac{l^{lk+1}-1}{l-1}\dim V$.
\begin{cor}
An algebraic supergroup $G$ is finite iff $G_{red}$ is finite or
iff $G_{ev}$ is finite.
\end{cor}
For any finite supergroup $G$ denote by $|G|$ the dimension of
$K[G]$. We call $|G|$ an {\it order} of $G$. By Lemma 6.2 of
\cite{w} there are pairwise-orthogonal idempotents $e_1,\ldots ,
e_n\in K[G]_0$ such that $\sum_{1\leq i\leq n}e_i=1$ and each
$K[G]e_i$ has a unique (nilpotent) maximal ideal $re_i, r=rad
K[G]$. Without loss of generality one can assume that
$\epsilon_G(e_1)=1$ and $\epsilon_G(e_i)=0$ for $i\geq 2$. A
supersubalgebra of $K[G]$, generated by $e_1,\ldots , e_n$, is
denoted by $B$.
\begin{lm}
Any idempotent of $K[G]$ belongs to $B$.
\end{lm}
Proof. Let $x=x_0+x_1$ be an idempotent in $K[G]$. The equality
$x^2=x$ implies that $x_0$ is also an idempotent and
$2x_0x_1=x_1$. Multiplying the last equality by $x_0$ we obtain
that $4x_0x_1=x_1$, hence $x_1=0$. Any idempotent $xe_i$ belongs
to $K[G]e_i$. In particular, $xe_i=\alpha e_i +ye_i$, where
$\alpha=0, 1$ and $y\in r$. On the other hand, $(x-\alpha)e_i$ is
again an idempotent that equals to a nilpotent element $ye_i$.
Hence any $xe_i$ is either zero or it is equal to $e_i$. Therefore
$x=\sum_{1\leq i\leq n} xe_i\in B$.
\begin{lm}
The algebra $B$ is a Hopf (super)subalgebra.
\end{lm}
Proof. We have to check that $B$ is a (super)subcoalgebra. The
radical of the superalgebra $A=K[G]\otimes K[G]$ equals to
$J=r\otimes K[N] +K[N]\otimes r$. Elements $e_i\otimes e_j$ are
pairwise orthogonal idempotents and their sum equals $1\otimes 1$.
Furthermore, $A e_i\otimes e_j /J e_i\otimes e_j\simeq K\otimes
K\simeq K$. Applying Lemma 7.6 to the finite supergroup $G\times
G$ we conclude that $\delta_G(B)\subseteq B\otimes B$.

Consider the natural epimorphism of supergroups $G\to SSp \ B$ and
denote by $G^{(0)}$ its kernel. We can consider $G^{(0)}$ as a
{\it connected component} of $G$. The equivalent definition will
be given in Section 9.

Assume that there is a finite supergroup $G$ and its supersubgroup
$H$ such that $\tilde{\tilde{G/H}}$ is not affine. The pair $(G,
H)$ is said {\it bad}. A bad pair defines a vector $v_{G, H}=(|G|,
|G|-|H|)\in {\bf N}^2$. Order vectors from ${\bf N}^2$
lexicographically from left to right. Choose a bad pair $(G, H)$
whose $v_{G, H}$ is minimal. Denote $K[G]^H$ by $R$.

By Theorem 5.2(ii) the property to be a faithfully exact
supersubgroup is transitive. Therefore, there is not any
supersubgroup $H'$ such that $H < H' < G$. On the other hand, the
superideal $K[G]R^+ + K[G]s_G(R^+)$ is contained in $I_H$ and it
is Hopf one. There is a supersubgroup $H'$ such that
$I_{H'}=K[G]R^+ + K[G]s_G(R^+)$. Since $H\leq H'$ it follows that
either $H=H'$, that is $I_H=K[G]R^+ + K[G]s_G(R^+)$, or $G=H'$,
that is $R=K$.

Assume that $E < G^{(0)} < G$. The minimality of $v_{G, H}$ and
Proposition 7.1 imply that $HG^{(0)} < G$ and therefore,
$H=HG^{(0)}$. Thus $G^{(0)}\leq H$. In the same way, by
Proposition 6.4 we see that $\tilde{\tilde{G/H}}$ is affine.
Finally, if $G^{(0)}=E$, then $G$ is pure even and
$\tilde{\tilde{G/H}}$ is always affine. The remaining case
$G=G^{(0)}$ means that $K[G]$ is a local superalgebra with
$\ker\epsilon_G=r$.
\begin{lm}
The superalgebra $K[G]$ is a free $R$-module.
\end{lm}
Proof. By Theorem 5.1 $K[G]$ is a flat$ R$-module. It remains to
notice that $R^+$ is a maximal nilpotent ideal of $R$ and use
Corollary 2.1 \cite{bur}, II, \S 3.
\begin{lm} Let $L$ be an affine supergroup and $N$ be an its
supersubgroup such that $I_N=K[L]T^+ + K[L]s_G(T^+), T=K[L]^N$.
Then, the canonical morphism $L\to SSp \ T$ induces an inclusion
$\tilde{\tilde{L/N}}\subseteq SSp \ T$.
\end{lm}
Proof. One has to check that the induced morphism $(L/N)_{(n)}\to
SSp \ T$ is injective. Let $g_1, g_2\in L(A), A\in SAlg_K$, and
assume that $g_1|_T=g_2|_T$. For a given $t\in T^+$ we have
$\delta_L(t)=t\otimes 1 +\sum h_1\otimes t_2$, where $t_2\in T^+$.
By definition,
$$(g_1^{-1}g_2)(t)=g_1(s_L(t))+\sum g_1(s_L(h_1))g_2(t_2)=
g_1(s_L(t))+\sum g_1(s_L(h_1))g_1(t_2)=$$
$$g_1(s_L(t)+\sum s_L(h_1)t_2)=g_1(\epsilon_L(t))=0.$$

Analogously, $$(g_1^-1 g_2)(s_L(t))=g_2(s_L(t))+\sum
(-1)^{|h_1||t_2|}g_1(t_2)g_2(s_L(h_1))=$$$$g_2(s_L(t))+\sum
(-1)^{|h_1||t_2|}g_2(t_2)g_2(s_L(h_1))=g_2^{-1}(\epsilon_L(t))=0.$$
It follows that $g_1^{-1}g_2(I_N)=0$, that is $g_1^{-1}g_2\in
N(A)$.

If $H=H'$, then Lemma 7.5 implies that the morphism of
$K$-functors $G\times H\to G\times_{SSp \ R} G$ defined as $(g,
h)\mapsto (g, gh), g\in G(A), h\in H(A), A\in SAlg_K,$ is an
isomorphism. Combining Proposition 4.2 with Lemma 7.4 we obtain
that $\tilde{\tilde{G/H}}\simeq SSp \ R$!

Let $R=K$. Since the ideal $I_H$ is nilpotent, one can repeat the
arguments from Propositions 6.5 and 6.6 to conclude that $G=H$.
Resuming all the above we obtain
\begin{ter}
Let $G$ be a finite supergroup and $H\leq G$. Then,
$\tilde{\tilde{G/H}}=SSp \ K[G]^H$.
\end{ter}
\begin{rem}
In the conditions of the above theorem we have
$\tilde{\tilde{G/H}}=\tilde{G/H}$. In fact, $K[G]^H$ is obviously
noetherian and it remains to use Proposition 4.2.
\end{rem}

\section{Brunden's question, $char K=p > 0$}

Let $G$ be an algebraic supergroup and $H$ be its supersubgroup
such that $H_{ev}$ is reductive. It is equivalent to say that
$H_{ev}$ is geometrically reductive \cite{hab, sesh} and the last
property is kept by taking quotients and normal subgroups
\cite{bs}. As above, we assume that $K$ is algebraically closed.
Let $A$ be a superalgebra. Denote by $A^{(n)}$ the superalgebra
that coincides as a ring with $A$, but where each $a\in K$ acts as
$a^{p^{-n}}$ does on $A$ (see \cite{jan}).
\begin{lm}(see \cite{mu}, 3.1(a))
If $A$ is a Hopf superalgebra and $char K=p >0$, then the linear
map $F^n : x\mapsto x^{p^n}$ (Frobenius morphism) is a Hopf
superalgebra morphism $A^{(n)}\to A$.
\end{lm}
Proof. The identity $a^{p^n}=a_0^{p^n}$ for $a\in A$ implies that
$F^n$ is a superalgebra morphism. Since $\delta_A$ and $s_A$ are
superalgebra morphisms, the equations $s_A F^n=F^n s_A$ and
$\delta_A F^n= (F^n\otimes F^n)\delta_A$ follow easily.

Denote by $f_n : G\to SSp \ K[G]^{(n)}$ the morphism of
supergroups dual to $F^n : K[G]^{(n)}\to K[G]$. The kernel
$G_n=\ker f_n$ is called the $n$-th {\it infinitesimal}
supersubgroup. By Theorem 6.1, $G_n$ is faithfully exact
supersubgroup of $G$. Besides, $\tilde{G/G_n}\simeq SSp
K[G]_0^{p^n}$. If $G$ is algebraic, then any $G_n$ is finite.
\begin{lm}
Let $L$ be an algebraic supergroup. For sufficiently large $t\geq
1$ the epimorphism $L\to \tilde{L/L_t}$ induces an epimorphism
$L_{ev}\to \tilde{L/L_t}$. In particular, if $L_{ev}$ is
reductive, then $\tilde{L/L_t}$ is also reductive.
\end{lm}
Proof. One has to check tat $K[L]_0^{p^t}\bigcap K[L]_1^2=0$ for
for some $t\geq 1$. Let $I$ be a radical of $K[L]_0$. The algebra
$K[L]_0$ is noetherian. It follows that $I^s=0$ for some $s\geq
1$. If $W$ is a complement of vector subspace $I$ to $K[L]_0$,
then $K[L]_0^{p^t}= W^{p^t}$, provided $p^t\geq s$. It remains to
notice that $K[L]_1^2\subseteq I$ and $W^{p^t}\bigcap I=0$.
\begin{pr}
The quotient $\tilde{\tilde{G/H}}$ is affine.
\end{pr}
Proof. Consider the supersubgroup $HG_t$. By Lemma 3.1 one can
assume that $H_{ev}\to \tilde{H/H_t}=\tilde{HG_t/G_t}$ is an
epimorphism. Combine with Corollary 4.5 from \cite{cps} we see
that the quotient of $\tilde{G/G_t}$ over the supersubgroup
$\tilde{HG_t/G_t}$ is affine. Since the property to be a
faithfully exact supersubgroup is transitive we refer to Theorem
7.1 and Proposition 7.1 to conclude the proof.

\section{Quotients by normal supersubgroups, $char K=0$}

An algebraic supergroup $G$ is called {\it pseudoconnected} if
$\bigcap_{n\geq 0}{\cal M}^n=0$, where ${\cal M}=\ker\epsilon_G$.
\begin{lm}
Let $G$ be an algebraic supergroup. The superideal
$I=\bigcap_{n\geq 0}{\cal M}^n$ is a Hopf superideal and the
supersubgroup $G^{[0]}=V(I)$ is normal and connected.
\end{lm}
Proof. By definition, $s_G({\cal M})={\cal M}$. It remains to
check that $I$ is a coideal and $\nu_l(I)\subseteq I\otimes K[G]$.
The trivial supersubgroup is obviously normal. In particular,
$\nu_l({\cal M})\subseteq {\cal M}\otimes K[G]$ which implies
$\nu_l({\cal M}^n)\subseteq {\cal M}^n\otimes K[G]$ for all $n\geq
0$ and we conclude that  $\nu_l(I)\subseteq I\otimes K[G]$.
Furthermore,
$$\delta_G({\cal M}^n)\subseteq \sum_{0\leq i\leq n}{\cal M}^i\otimes
{\cal M}^{n-i}\subseteq\bigcap_{0\leq i\leq n}({\cal M}^i\otimes
K[G] + K[G]\otimes {\cal M}^{n-i})$$ and
$$\delta_G(I)\subseteq\bigcap_{n\geq 0}\delta_G({\cal M}^n)\subseteq
\bigcap_{n\geq 0}({\cal M}^n\otimes K[G] + K[G]\otimes {\cal
M}^n)= I\otimes K[G] +K[G]\otimes I .$$

The supersubgroup $G^{[0]}$ is called a {\it pseudoconnected
component} of $G$. It is clear that $G$ is pseudoconnected iff
$G=G^{[0]}$. Additionally, $Dist(G)=Dist(G^{[0]})$ and
$Lie(G)=Lie(G^{[0]})$. The proof of Lemma 9.1 shows that if
$N\unlhd G$, then $N^{[0]}\unlhd G$ also. Besides, an epimorphic
image of a pseudoconnected supergroup is again pseudoconnected.
\begin{lm}(Krull's intersection theorem)
Let $A$ be a finitely generated commutative superalgebra and $V$
be a finitely generated $A$-supermodule. For any superideal $I$ of
$A$ we have $\bigcap_{t\geq 0}I^t V=\{v\in V|\mbox{there exists}
 \ x\in I_0 \ \mbox{such that} \ (1-x)v=0\}$.
\end{lm}
Proof. Observe that $V$ is finitely generated as a $A_0$-module.
Since $I^t\subseteq I_0^{[\frac{t}{2}]}\bigoplus
I_0^{[\frac{t-1}{2}]}I_1$ we see that $I^t V\subseteq
I_0^{[\frac{t-1}{2}]}V_0\bigoplus I_0^{[\frac{t-1}{2}]}V_1$.
Proposition 5, \cite{bur}, III, \S3, concludes the proof.
\begin{pr}
Let $\pi : G\to H$ is an epimorphism of algebraic supergroups. If
$char K=0$, then the induced short sequence of Lie superalgebras
$$0\to Lie(\ker\pi)\to Lie(G)\stackrel{d\pi}{\to} Lie(H)\to 0$$
is exact.
\end{pr}
Proof. Since $I_{\ker\pi}=K[G]\pi^*(\ker\epsilon_H)$ it obviously
implies that $\ker d\pi=Lie(\ker\pi)$. Combining Lemma 9.2 with
word-by-word repeating the proof of Proposition 7.6, \cite{jan},
part I, we obtain that $d\pi : Dist(G)\to Dist(H)$ is surjective.
Now Lemma 3.1 concludes the proof.
\begin{rem}
If $char K=0$, then any algebraic supergroup $G$ is reduced (or
smooth), meaning that the radical of $K[G]$ coincides with
$I_{G_{ev}}$. Indeed, by Theorem 11.4 of \cite{w} the Hopf
superalgebra $K[G]/I_{G_{ev}}$ is reduced, and in particular,
$G_r=G_{ev}$. It is proved in \cite{f} that a completion of $K[G]$
with respect to a ${\cal M}$-adic topology is isomorphic to
$K[[t_1,\ldots , t_m|z_1 ,\ldots z_n]]$ for any ${\cal M}\in
Max(K[G])$, where $m=\dim Lie(G)_0, n=\dim Lie(G)_1$ (see also
Lemma 3.1 and use automorphism $T_{\cal M}$ from Theorem 5.1).
\end{rem}
Any (left) $G$-supermodule $V$ is a $Dist(G)$-supermodule via
$\phi v=\sum (-1)^{|\phi||v_1|}v_1\phi(f_2)$ for $\phi\in Dist(G)$
and $v\in V, \tau_V(v)=\sum v_1\otimes f_2$, see \cite{bk, bkl}.
If $V, V'$ are $G$-supermodules, then we have a canonical
embedding $Hom_G(V, V')\subseteq Hom_{Dist(G)}(V, V')$.

The proof of the following lemmas is a copy of the proofs of
Proposition 7.5, Lemmas 7.15 and 7.16, \cite{jan}, part I.
\begin{lm}
Let $H_1, H_2$ be supersubgroups of an algebraic supergroup $G$
and $H_1$ be pseudoconnected. Then $H_1\subseteq H_2$ is
equivalent to $Dist(H_1)\subseteq Dist(H_2)$. Additionally, if
$char K=0$, then $H_1\subseteq H_2$ is equivalent to
$Lie(H_1)\subseteq Lie(H_2)$.
\end{lm}
\begin{lm}
If $G$ is pseudoconnected and algebraic, then a supersubspace $W$
of a $G$-supermodule $V$ is a $G$-supesubmodule iff $W$ is a
$Dist(G)$-supersubmodule. If $char K=0$, then $W$ is a
$G$-supesubmodule iff $W$ is a $Lie(G)$-supersubmodule.
\end{lm}
\begin{lm}
If $G$ is pseudoconnected and algebraic, then $Hom_G(V,
V')=Hom_{Dist(G)}(V, V')$ for any $G$-supermodules $V$ and $V'$.
If $char K=0$, then $Hom_G(V, V')=Hom_{Lie(G)}(V, V')$.
\end{lm}
\begin{pr}
If $L$ is algebraic and $N\unlhd L$, then the equality
$Lie(N)=Lie(L)$ implies that $\tilde{\tilde{L/N}}=\tilde{L/N}$ is
affine and finite.
\end{pr}
Proof.  Without loss of generality, one can assume that $K[L]^N=K$
and $K$ is algebraically closed. Denote the supergroup $N\bigcap
G_{red}$ by $U$. We have
$$Lie(U)=Lie(N)\bigcap Lie(G_{red})=Lie(G)\bigcap
Lie(G_{red})=Lie(G_{red}).$$ It follows that the affine group
$\tilde{L_r/U}$ is finite. In other words, an algebra
$B=K[\tilde{L_{red}/U}]=K[L_{red}]^U$ is finite dimensional. By
Lemma 6.2 of \cite{w} we have $B=\prod_{1\leq i\leq n} Be'_i$,
where $e'_1,\ldots e'_n$ are pairwise orthogonal idempotents such
that $\sum_{1\leq i\leq n}e'_i=1$. Besides, each algebra $Be'_i$
is isomorphic to $K$. By Corollary 1, \cite{bur}, II, \S 4, there
are pairwise orthogonal idempotents  $e_1,\ldots , e_n\in K[L]_0$
such that their respective images in $K[L_{red}]$ coincide with
$e'_1,\ldots e'_n$ and $\sum_{1\leq i\leq n}e_i=1$. Consider an
idempotent $e=e_i$. As in Lemma  6.3 we have $\delta_L(e)-e\otimes
1\in r\otimes K[L] +K[L]\otimes (r+I_N)$, where $r=rad K[L]$. On
the other hand, for any odd exponent $k$, the equality
$(\delta_L(e)-e\otimes 1)^k=\delta_L(e)-e\otimes 1$ holds. For
sufficiently large (odd) integer $k$ we infer
$\delta_L(e)-e\otimes 1\in K[L]\otimes I_N$, forcing $e\in K$. It
obviously implies that $n=1, e_1=1$ and $B=K$. Repeating the
arguments from Propositions 6.5 and 6.6, we see that
$\tilde{\tilde{L/N}}$ is affine. By Remark 6.1 $R=K[L]^N$ is
finitely generated and $\tilde{\tilde{L/N}}=\tilde{L/N}$. Let $I$
be a radical of $R$. Since $K[L]$ is a faithfully flat $R$-module
we have $r\bigcap R\subseteq I=R\bigcap K[L]I\subseteq r ,$ that
is $I=R\bigcap r$. In other words, the induced morphism
$L_{red}\to (\tilde{L/N})_{red}$ is an epimorphism. In particular,
$(\tilde{L/N})_{red}\simeq\tilde{L_{red}/U}$ is finite.
\begin{cor}
For any algebraic supergroup $L$ the quotient $\tilde{L/L^{[0]}}$
is a finite supergroup.
\end{cor}
\begin{rem}
If $char K=0$ and $Lie(G_{ev})=Lie(N_{ev})$, then all statements
of Proposition 9.2 hold also.
\end{rem}
We define a {\it connected} component of an algebraic group $G$ as
the preimage of $(\tilde{G/G^{[0]}})^{(0)}$ in $G$.  It is not
difficult to check that $G^{(0)}$ can be defined as the
intersection of kernels of all morphisms $G\to L$, where $L$ is an
even etale (super)group. Notice also that for any finite $G$ both
components are the same. In fact, if $e_1,\ldots , e_n$ are all
primitive idempotents of $K[G]$ and $\epsilon_G(e_1)=1$, then
$\ker\epsilon_G=r + K[G]e$, where $e=\sum_{2\leq i\leq n} e_i$ and
$r=rad K[G]$. Thus $I_{G^{[0]}}=K[G]e=I_{G^{(0)}}$.
\begin{rem}
The supersubgrop $G^{(0)}$ is an open subfunctor of $G$. In fact,
if $e_1,\ldots , e_n$ are all primitive idempotents of
$K[G]^{G^{[0]}}$, then $G^{(0)}=G_{e_1}$, provided
$\epsilon_G(e_1)=1$.
\end{rem}
\begin{quest}
Does $G^{[0]}$ coincide with $G^{(0)}$ for arbitrary $G$?
\end{quest}
\begin{quest}
Is $G^{[0]}$ an open subfunctor for arbitrary $G$?
\end{quest}
\begin{lm}
If $G$ is pseudoconnected or connected, then $Lie(G)=0$ implies
$G=E$. In particular, if $char K=0$ and $G$ is algebraic, then
$G^{(0)}=G^{[0]}$.
\end{lm}
Proof. In the above notations $Lie(G)=0$ iff $r\subseteq K[G]e$,
that is $\ker\epsilon_G=Be$ is an idempotent and nilpotent ideal
simultaneously. As for the second statement, Proposition 9.1
implies $Lie(G^{(0)}/G^{[0]})=0$.
\begin{lm}
If $H$ is a supersubgroup of $G$, then
$Lie(N_G(H))=(Lie(G)/Lie(H))^{{\bf Ad}(H)}$.
\end{lm}
Proof. As in \cite{dg}, II, \S 5, Lemma 5.7, it is enough to
observe that $x\in {\bf Lie}(N_G(H))(K)$ iff for any $A\in SAlg_K$
and $h\in H(A)$ we have
$$e^{\varepsilon_{|x|}x'} h' e^{-\varepsilon_{|x|}x'} h'^{-1}=e^{\varepsilon_{|x|}(x'-{\bf Ad}(h)(x'))}
\in H(A[\varepsilon_0,\varepsilon_1]),$$
$$e^{-\varepsilon_{|x|}x'} h' e^{\varepsilon_{|x|}x'}
h'^{-1}=e^{\varepsilon_{|x|}({\bf Ad}(h)(x')-x')} \in
H(A[\varepsilon_0,\varepsilon_1]),$$ where $x'=G(\iota^A_K)(x),
h'= G(p_A)(h)$. In other words, $x\in {\bf Lie}(N_G(H))(K)$ iff
$x'\pm {\bf Ad}(h)(x')\in {\bf Lie}(H)(A)$ for any superalgebra
$A$.

Let us return to the situation of Section 6, that is $N\leq G$,
where $G$ is algebraic. As before, one can assume that $K$ is an
algebraically closed of zero characteristic and, if it is
necessary, that $K[G]^N=K$. Define the {\it lower central} ({\it
solvable}) series of $L$ by $L^1=L, L^{i+1}=[L^i , L]$
(respectively, $L^{(0)}=L, L^{(i+1)}=[L^{(i)} , L^{(i)}]$.

As in Section 7 a pair $(G, N)$ is called {\it bad}, whenever
$\tilde{\tilde{G/N}}$ is not affine. A vector $v\in {\bf N}^2$ is
called {\it positive} iff at least one coordinate of $v$ is
positive. Partially order the set of bad pairs by $(G, N) < (G' ,
N')$ iff $\mbox{s}\dim Lie(G')-\mbox{s}\dim Lie(G)$ is positive,
otherwise  $\mbox{s}\dim Lie(G')-\mbox{s}\dim Lie(G)=0$ and
$\mbox{s}\dim Lie(N')-\mbox{s}\dim Lie(N)$ is positive. Choose a
minimal bad pair $(G, N)$.
\begin{pr}
If the superalgebra $L=Lie(G)$ is not simple, then it is either
semisimple, whose a unique proper ideal is $L^2$, or $L^2=0$.
\end{pr}
Proof. Let $I$ be an proper ideal of $L$. Consider $L$ as
$G$-supermodule via ${\bf Ad} : G\to GL(L)$. By Lemma 9.4 $I$ is a
$G$-supersubmodule and we define the induced morphism ${\bf Ad}_I
: G\to GL(L/I)$. Denote $\ker {\bf Ad}_I$ by $H$. By Proposition
9.1 $Lie(H)$ is a proper supersubalgebra of $L$ iff
$L^2\not\subseteq I$. As above, the minimality of $(G, N)$ and
Propositions 6.4 and 7.1 imply that $L^2$ is a smallest (possibly
zero) ideal of $L$. If $L$ is not semisimple, then considering the
morphism $Ad|_I : G\to GL(I)$ for a proper abelian ideal $I$, we
obtain $[L, I]=0$. In particular, $L^3 =0$. Finally, repeating the
above arguments for the morphism ${\bf Ad}$ we conclude that
$L^2=0$.

Using Propositions 6.4 and 7.1 and Theorem 7.1 as well, one can
always assume that $G$ is connected. If $L^2=0$, then $Dist(G)$ is
a commutative superalgebra. In its turn, $K[G]$ is cocommutative
and $G$ is an abelian supergroup. In particular, $G_{ev}\unlhd G$
and the minimality arguments imply that $Lie(G_{ev})=0$. Corollary
7.1 and Theorem 7.1 conclude the proof in this case. If $L$ is
simple, then $Lie(N)=0$. By Proposition 6.4 and Lemma 9.6 this
case is reduced to $N=N_{ev}$. The algebra $K[G]_0$ can be
regarded as a coordinate algebra of an affine scheme on which $N$
acts on the right. Combining  \cite{dg}, III, \S2,
$\mbox{n}^{\circ}$4, with $K[G]^N=K$ we obtain that $K[G]_0$ is
finitely dimensional. By Theorem 7.1 this case is also done.

Finally, if $L$ is semisimple, then we consider the supersubgroup
$H=(G_{ev}N)^{[0]}$. Denote by $Q$ the normalizer $N_G(H)$ and by
$D$ its Lie superalgebra. By the above, one can assume that
$Lie(N)\neq 0$. It follows that $L^2\subseteq Lie(N)$. Denote
$Lie(G_{ev}N)=L_0 + Lie(N)$ by $M$. Since $L^2\subseteq M$ the
supergroup $H$ acts trivially on $L/M$. By Lemma 9.7
$D/M=(L/M)^{{\bf Ad}(H)}=L/M$ that infers $G=Q$ and $H\unlhd G$.
Again, by the minimality of $(G, N)$ we have either $G=G_{ev}N$
and then
$\tilde{\tilde{G/N}}\simeq\tilde{\tilde{G_{ev}/(G_{ev}\bigcap
N)}}$ is affine, or $L_0\subseteq Lie(N)$ and Remark 9.2 concludes
the proof of Theorem 6.2.

\section{Two examples}

One more example of not necessary normal but faithfully exact
supersubgroup is given by a {\it Levi supersubgroup}. In notations
of Proposition 6.3, a Levi supersubgroup $L_{s,t}$ of $G=GL(m|n)$
consists of all matrices satisfying the equations $A_{12}=0,
A_{21}=0, B_{12}=0, B_{21}=0, B_{22}=0, C_{12}=0, C_{21}=0,
C_{22}=0, D_{12}=0, D_{21}=0$. Additionally,  the blocks $A_{22}$
and $D_{22}$ are diagonal matrices. It is clear that
$L_{s,t}\simeq GL(s|t)\times T$, where $T$ is a torus of dimension
$m+n-s-t$. Represent the coordinate superalgebra of $K[GL(m|n)]$
as $K[A, B, C, D]_{d_1 d_2}$, where the blocks $A, B, C, D$ are
identified with the sets of their coefficients, and $d_1=\det(A),
d_2=\det(D)$.

Define the map $\pi : \underline{m+n}\to \underline{m+n}$ by the
rule $\pi(i)=i+n$ and $\pi(j)=j-m$ for $1\leq i\leq m < j\leq
m+n$.
\begin{lm}
There is a canonical isomorphism $\psi : GL(m|n)\to GL(n|m)$ such
that $\psi(L_{s, t})=L_{t, s}$.
\end{lm}
Proof. Denote the matrix coordinate functions on $GL(m|n)$ by
$a_{ij}$ and the similar functions on $GL(n|m)$ by $a'_{ij}$. It
is easy to see that $\psi^*(a_{ij})= a'_{\pi(i), \ \pi(j)}$
induces the required Hopf superalgebra isomorphism
$K[GL(m|n)]\simeq K[GL(n|m)]$. In fact, only the equality $\psi^*
s_{GL(m|n)}= s_{GL(n|m)}\psi^*$ is not trivial. It is enough to
prove it for generators $a_{ij}$ and using the following formulaes
$$\psi^*\left(\begin{array}{cc}
A & B \\
C & D \end{array}\right)=\left(\begin{array}{cc}
\psi^*(D) & \psi^*(C) \\
\psi^*(B) & \psi^*(A) \end{array}\right),$$
$$s_{GL(m|n)}\left(\begin{array}{cc}
A & B \\
C & D \end{array}\right)=\left(\begin{array}{cc}
(A-BD^{-1}C)^{-1} & -A^{-1}B(D-CA^{-1}B)^{-1} \\
-D^{-1}C(A-BD^{-1}C)^{-1} & (D-CA^{-1}B)^{-1} \end{array}\right)
$$
it can be done by straightforward calculations.

Using Lemma 10.1 and transitivity of inducing functor it remains
to prove that $L_{m, \ n-1}$ (or $L_{m-1, \ n}$) is faithfully
exact. In what follows let $G=GL(m|n)$.
\begin{lm}
The elements $a_{i, \ m+n}s_G(a_{m+n, \ j})$ for $1\leq i, j\leq
m+n,$ are $L=L_{m, \ n-1}$-invariants. Moreover, they generate a
left coideal supersubalgebra $R$ such that $I_L=K[G]R^+$.
\end{lm}
Proof. Denote the image of $x\in K[G]$ in $K[L]$ by
$\overline{x}$. Since
$$\delta_G(a_{i, \ m+n}s_G(a_{m+n, \ j}))=
\sum_{1\leq t, \ l\leq m+n} (-1)^{|a_{lj}||a_{m+n, \ l}|+|a_{t, \
m+n}||a_{lj}|}a_{it}s_G(a_{lj})\otimes a_{t, \ m+n}s_G(a_{m+n, \
l}),$$ it follows that $R$ is a left coideal. Considering $K[G]$
as a $K[L]$-supercomodule (where $L$ acts on $G$ by the right
multiplication) we determine
$$\tau_{K[G]}(a_{i, \ m+n}s_G(a_{m+n, \ j}))=a_{i, \ m+n}s_G(a_{m+n, \ j})\otimes
\overline{a_{m+n, \ m+n}s_G(a_{m+n, \ m+n})}=$$
$$a_{i, \ m+n}s_G(a_{m+n, \ j})\otimes 1 .$$

Notice that $a_{i, \ m+n}s_G(a_{m+n, \ j})-\delta_{i, \
m+n}\delta_{m+n, \ j}\in I_L$ for $1\leq i, j\leq m+n$ and
consider an index $i\neq m+n$. We have $a_{i, \ m+n}=\sum_{1\leq
j\leq m+n} a_{i, \ m+n}s_G(a_{m+n, \ j})a_{j, \ m+n}\in I_L$ and
symmetrically, $s_G(a_{m+n, \ j})\in I_L$ for $j\neq m+n$. Denote
by $I$ a superideal generated by elements $a_{i, \ m+n},
s_G(a_{m+n, \ j})$ for $1\leq i, j < m+n$. It is obvious that
$L\subseteq V(I)$. On the other hand, $g\in V(I)(A)$ iff $g\in
Stab_G(Kv_{m+n})(A)$ and $g^{-1}\in Stab_G(\sum_{1\leq i\leq
m+n-1} Kv_i)(A)$. The superversion of \cite{jan}, part I (1.4)
completes the proof of this claim.

By Proposition 5.2 and Theorem 5.1 all we have to show is that
$K[G]{\cal M}\neq K[G]$ for all ${\cal M}\in Max(R)$. Using the
reduction from Corollary 5.1 one can work with algebraic groups
$L_{0, \ n-1}\subseteq GL(0|n)\simeq GL(n)$. In other words, we
can set $G=GL(n), \ L=Stab_G (Kv_n)\bigcap Stab_G(\sum_{1\leq
i\leq n-1}Kv_i)$. By Corollary 4.5 of \cite{cps}, the quotient
$\tilde{G/L}$ is affine and it is isomorphic to $Sp_K \ K[G]^L$.
In particular, $K[G]^L\leq K[G]$. We will show that $K[G]^L$ is
generated by the elements $a_{in}s_G(a_{nj})$ for $1\leq i, j\leq
n$.

Consider an element $\frac{f}{d^k}\in K[G]$, where $f\in K[a_{ij}|
1\leq i, j\leq n]$ and $d=\det(a_{ij})$. Represent $f$ as
$$f=\sum_{\alpha\in {\bf N}^n } f_{\alpha}\prod_{1\leq i\leq
n}a_{in}^{\alpha_i},\mbox{ where } f_{\alpha}\in K[a_{ij} | 1\leq
i\leq n, \ 1\leq j\leq n-1].$$

It can be checked easily that $\frac{f}{d^k}\in K[G]^L$ iff all
monomials $\prod_{1\leq i\leq n}a_{in}^{\alpha_i}$ in the above
representation of $f$ have degree $k$ and  all coefficients
$f_{\alpha}$ are $GL(n-1)$-semi-invariants of weight $k$ with
respect to the action of this group by right multiplications on
the variety of $n\times (n-1)$-matrices $M_{n\times (n-1)}$. By
Igusa's Theorem (see \cite{ig} and Corollary 3.5 of \cite{cep})
the algebra $K[M_{n\times (n-1)}]^{SL(n-1)}=K[a_{ij}| 1\leq i\leq
n, \ 1\leq j\leq n-1]^{SL(n-1)}$ is generated by all minors of
size $n-1$ which are just semi-invariants $s_G(a_{nj})d$ of weight
1.

The next example shows that a quotient $\tilde{\tilde{X/G}}$ (as
well $\tilde{X/G}$) is not always affine even if $G$ is finite. In
what follows $char K=0$. Let $G=G_a^{-}$ be an odd unipotent
supergroup such that $K[G]=K[t]$, where $|t|=1$,
$\delta_G(t)=t\otimes 1 +1\otimes t$ and $\epsilon_G(t)=0,
s_G(t)=-t$. A superspace $V$ is a $G$-supermodule iff there is
$\phi\in End_K(V)_1, \phi^2=0,$ such that $\tau_V(v)=v\otimes 1
+\phi(v)\otimes t$. Moreover, $V^G=\ker\phi$. Assume that
$\mbox{s}\dim V=(1, 1)$ and $v_1, v_2$ form a ${\bf
Z}_2$-homogeneous basis of $V$. Set $\phi(v_1)=v_2, \phi(v_2)=0$.
The symmetric superalgebra $S(V)$ has the induced $G$-supermodule
structure by
$$\tau_{S(V)}(v_1^r)=v_1^r\otimes 1 +rv_1^{r-1}v_2\otimes t, \
\tau_{S(V)}(v_1^{r-1}v_2)=v_1^{r-1}v_2\otimes 1, r\geq 0.$$ Since
$\tau_{S(V)}$ is a superalgebra morphism, there is an affine
scheme $X$ such that $K[X]=S(V)$ and $G$ acts on $X$. Denote
$K[X]^G$ by $R$. We have $R_0=K, R_1=\sum_{t\geq 0}Kv_1^t v_2$ and
$R_1^2=0$. In particular,  $R$ is commutative as an algebra.
\begin{lm}
Any flat $R$-supermodule is flat as a module.
\end{lm}
Proof. Any exact sequence of $R$-modules $0\to V\to W\to U\to 0$
can be turn into an exact sequence of $R$-supermodules. In fact,
set $V_1=V\bigcap R_1W$ and let $V_0$ be a complement of $V_1$ to
$V$. Since $V_0\bigcap R_1W=0$, then $R_1W$ has a complement $W_0$
to $W$ such that $V_0\subseteq W_0$. Finally, $U= W_0/V_0\bigoplus
W_1/V_1$. If $M$ is a flat $R$-supermodule, then the functor
$M\otimes_R ?$ takes our sequence (of supermodules!) to an exact
one.

Assume that $\tilde{X/G}$ or $\tilde{\tilde{X/G}}$  is an affine
superscheme. By Proposition 4.1 $S(V)$ is a faithfully flat
$R$-supermodule. Combining Lemma 10.3  with Corollary 2.1
\cite{bur}, II, \S 3, we obtain that $S(V)$ is a free $R$-module.
Let $Rf\simeq R, f\in S(V)$. The equality $Rf=Rf_0 + Kf_1$ implies
that either $Rf=Rf_0$ or $Kf_1$ is a projective $R$-module. Again
by Corollary 2.1 \cite{bur}, II, \S 3, $Kf_1$ has to be free, but
it is  obviously impossible. Thus $S(V)=\bigoplus Rf$, where $f$
runs over a basis of $S(V)_0$. Without loss of generality one can
assume that some $f$ equals $1$. On the other hand, the summand
$R$ has nontrivial intersections with all other summands! This
contradiction shows that both $\tilde{X/G}$ and
$\tilde{\tilde{X/G}}$ are not affine.

\begin{center}
\bf Acknowledgements
\end{center}
The author thanks FAPESP (proc. 07/54834-9) for the support during
his visit to Sao Paulo University. This work was also partially
supported by RFFI 07-01-00392. Special thanks to Professor
Alexandr Grishkov for his invitation and to Professor Jonathan
Brundan for bringing to my attention the example of Levi
supersubgroup and communicating to me interesting questions.

\end{document}